\documentclass[10pt,a4paper]{article}
\usepackage{amsmath}
\usepackage{amssymb}
\usepackage{amsfonts}

\usepackage{geometry}
\usepackage{indentfirst}

\usepackage{tabularx}
\usepackage{booktabs}
\usepackage{graphicx}

\usepackage{graphicx}

\newtheorem{remark}{Remark}

\newcommand{\RR}{\mathbb{R}}

\newcommand{\dt}{\Delta t}
\def\epsi{{\varepsilon}}
\def\u{{\rho}}
\def\v{{j}}

\newcommand{\emat}{\end{pmatrix}}

\begin{document}

\begin{center}

  \LARGE
Asymptotic-Preserving Monte Carlo methods for transport equations in the diffusive limit
  \bigskip
  \medskip

  \large
  G. Dimarco\textsuperscript{1}, L. Pareschi\textsuperscript{1} and G. Samaey\textsuperscript{2}

  \medskip
  \smallskip

  \normalsize
  \textsuperscript{1}Department of Mathematics and Computer Science, University of Ferrara, Italy \\
  \textsuperscript{2}NUMA (Numerical Analysis and Applied Mathematics), Department of Computer Science, KU Leuven,
Celestijnenlaan 200A, 3001 Leuven, Belgium\\
\end{center}


%
%

\begin{abstract}	
We develop a new Monte Carlo method that solves hyperbolic transport equations with stiff terms, characterized by a (small) scaling parameter. In particular, we focus on systems which lead to a reduced problem of parabolic type in the limit when the scaling parameter tends to zero. Classical Monte Carlo methods suffer of severe time step limitations in these situations, due to the fact that the characteristic speeds 
go to infinity in the diffusion limit. This makes the problem a real challenge, since the scaling parameter may differ by
several orders of magnitude in the domain. To circumvent these time step limitations, we construct a new, asymptotic-preserving Monte Carlo method
that is stable independently of the scaling parameter and degenerates to a standard probabilistic approach for solving the limiting equation in the diffusion limit. The method uses an implicit time discretization to formulate a modified equation in which the characteristic speeds do not grow indefinitely when the scaling factor tends to zero. The resulting modified equation can readily be discretized by a Monte Carlo scheme, in which the particles combine a finite propagation speed with a time-step dependent diffusion term. We show the performance of the method
 by comparing it with standard (deterministic) approaches in the literature.
\end{abstract}

ving schemes.
 
 

\section{Introduction}
In many situations in which hyperbolic transport equations intervene, such as neutron transport \cite{Ce}, radiative transfer \cite{bird, Case, Cha}, plasma physics \cite{birsdall} or semiconductor simulation \cite{Mark}, the multiscale nature of the phenomena involved causes large difficulties for the development of efficient numerical methods. In fact, the scaling parameters that characterize the relevant time scales which determine evolution of such problems may differ by several
order of magnitude, making the problem very stiff \cite{Ce}. This stiffness limits the maximal time-step that can be taken when using an explicit discretization. In addition, in the limit when the scaling parameter tends to zero, the equations that are used to model these phenomena may change character,
passing from a hyperbolic to a parabolic structure \cite{BGL}. In the particular setting of a diffusive scaling \cite{BGL, CIP} that is considered in this work, the characteristic speeds of the hyperbolic system grow to infinity as the scaling parameter tends to zero, causing severe CFL restrictions in standard explicit numerical methods, both in  deterministic (grid-based) and stochastic Monte Carlo (particle-based) approaches. 

Thus, the first idea to deal with such systems consists in using implicit time integration methods. Unfortunately,  these techniques are only possible to implement for grid-based space discretizations. Moreover, the main drawback of such approaches is that they often lead to large nonlinear systems very hard to invert due to the high dimension of the equations involved. For this reason, it is desirable to develop numerical methods that are not fully implicit, yet able to overcome the computational cost caused by the stiffness by avoiding time step limitations related to the scaling parameter. A particularly appealing class of schemes, developed for facing these situations, are the so-called asymptotic-preserving schemes \cite{BEN, BPR12, BC, CL, DP1, GOS1, J, JPT1, JPT2, Klar, Kl2, Larsen, LemMieu2008, NP2}, which degenerate to consistent discretization of the limiting problem when the scaling parameter is set to zero. Recently, this approach has been considered in the framework of Implicit-Explicit (IMEX) Runge-Kutta schemes with the aim of deriving
high order numerical methods that are accurate in all regimes, i.e., regardless of the value of the scaling parameter  \cite{DP1, BPR12, BPR17, Risp}. In particular, these schemes also preserve the desired order of accuracy, even in the limit when the scaling parameter tends to zero. Recently, via projective integration \cite{GearKevrekidis03}, also fully explicit methods for stiff hyperbolic transport equations have been developed \cite{LafitteSamaey,LafitteLejonSamaey,LafitteMelisSamaey}.

In this work, we concentrate on the development of asymptotic-preserving Monte Carlo methods for solving hyperbolic transport equations in the diffusion limit. The Monte Carlo approach represents a very popular method to deal with transport equations due to its flexibility and low computational cost \cite{Cf, bird, Nanbu80, PRMC, sliu}. However, one of the main limitation of Monte Carlo methods is the difficulty to construct schemes which work uniformly for all values of the scaling parameters \cite{Cf,PRMC}. For this reason, various modified Monte Carlo techniques have been recently proposed in the case of the so-called hydrodynamic scaling \cite{dimarco5, dimarco6, dimarco1, CPmc, PRMC, Trazzi}. To the best of our knowledge, up to now, a Monte Carlo scheme which is able to work uniformly in the diffusive scaling without causing severe time step limitations has not been constructed. The main motivation for not using particle methods in this framework is that the characteristic speeds grow to infinity in the diffusive limit. 
Thus, the dominant technique to overcome this problem consists nowadays in employing domain decomposition strategies, in which a Monte Carlo discretization of the hyperbolic transport equation is solved in regions in which the scaling parameter is large, and the limiting equation is solved (deterministically) in regions where the scaling parameter is small. These approaches have been largely studied for kinetic equations both for the 
diffusive \cite{BAL,DEG1,Klar} and for the hydrodynamic scaling \cite{DIM1, DIM2}; see also, e.g., \cite{Garcia-AMAR} for related ideas. However, even if these methods are very efficient, they are affected by some difficulties due to the fact that it is not always a simple task to define the different regions of the domain in which the use of a macroscopic model is fully justified. 

In this paper, our main goal is to develop an asymptotic-preserving Monte Carlo method in the diffusive scaling. The method does not rely on domain decomposition strategies, nor on coupling with a deterministic discretization of the limiting equation. Instead, to overcome the parabolic stiffness, the main ingredients are a suitable reformulation of the original system based on an implicit time discretization \cite{BPR12,BPR17}, which leads to a modified equation where the characteristic speeds are bounded in terms of the scaling parameters. An appropriate splitting strategy for this equation then permits the construction of a Monte Carlo scheme that works independently of the value of the scaling parameter and that automatically degenerates to a classical random walk method for limiting the diffusion equation. The method is first constructed by using the Goldstein--€"Taylor model \cite{Goldstein} under the diffusive scaling. Successively, the scheme is extended to the case of the kinetic radiative transport equation \cite{Ce}.

The remainder of this paper is organized as follows. In Section~\ref{G-T}, we discuss the development of the new Monte Carlo method in the case of the Goldstein-Taylor model. First we introduce the model and its diffusion limit. We also recall the standard splitting method that leads to classical Monte Carlo schemes for kinetic equations. Next, we introduce the new asymptotic preserving Monte Carlo method. Subsequently, in Section~\ref{R-T}, we extend our methodology to the case of the radiative transport equation. After the introduction of the kinetic model in the diffusive scaling we again discuss its diffusive limit and the corresponding classical Monte Carlo schemes. Then, we then extend the asymptotic preserving Monte Carlo scheme to the case of the radiative transport. We present several numerical tests and analyze the performance of the new methods in Section~\ref{Num}. Finally, in Section~\ref{Concl}, we draw some conclusions and give an outlook to future research directions.

\section{Asymptotic-preserving Monte Carlo methods for the Goldstein Taylor model}
\label{G-T}
In this section we discuss the construction of asymptotic-preserving Monte Carlo methods using as a prototype problem the Goldstein-Taylor (GT) model in the diffusive scaling. To this aim, we first introduce the prototype problem and emphasize the drawbacks of a standard Monte Carlo approach. Next, by means of a suitable problem reformulation we construct our novel class of asymptotic preserving Monte Carlo (APMC) methods.

\subsection{The GT model in the diffusive limit\label{sec:model}}



The GT model \cite{Goldstein} can be written as the following system of kinetic equations 
\begin{equation}
 \left\{  
\begin{array}{l} 
\displaystyle  
\partial_t f_+ + {\mathcal V}\partial_x f_+ =\frac{\sigma}{2}\left(f_--f_+\right), \\[+.25cm]
\displaystyle    
\partial_t f_- - {\mathcal V}\partial_x f_- =\frac{\sigma}{2}\left(f_+-f_-\right),
\end{array}
\right. 
\label{I62a}
\end{equation}
where $f_{\pm}=f_{\pm}(x,t)\geq 0$, $x\in\Omega\subset \mathbb{R}$ and $t \in \mathbb{R}^+$.
The model describes the evolution of two sets of particles: particles with (constant) velocity ${\mathcal V}$, with density $f_+$, and particles with (constant) velocity $-{\mathcal V}$ with density $f_-$. Equation~\eqref{I62a} gives a physically intuitive description of the process: the left hand side denotes transport with the characteristic speeds $\pm\mathcal{V}$, whereas the right hand side encodes random velocity changes that can be interpreted as collisions: in both populations particles disappear with a rate $\sigma$ and reappear with the opposite velocity. 

The diffusive scaling corresponds to take ${\mathcal V}=1/\varepsilon$ and $\sigma=1/\varepsilon^2$ where $\varepsilon>0$ is the scaling parameter. The scaled model reads
\begin{equation}
 \left\{  
\begin{array}{l} 
\displaystyle  
\partial_t f_+ + \frac{1}{\varepsilon}\partial_x f_+ =\frac{1}{\varepsilon^2}\left(\dfrac{\u}{2}-f_+\right), \\[+.25cm]
\displaystyle    
\partial_t f_- - \frac{1}{\varepsilon}\partial_x f_- =\frac{1}{\varepsilon^2}\left(\dfrac{\u}{2}-f_-\right),
\end{array}
\right. 
\label{I62}
\end{equation}
where we have introduced the mass density $\rho(x,t)$ and the scaled momentum $j(x,t)$ defined as
\begin{equation}\label{eq:trasform}
\u = f_+ + f_-, \qquad \v = \dfrac{f_+ - f_-}{\varepsilon}.
\end{equation}
Now, using the above macroscopic variables system~\eqref{I62} can be written equivalently in the form
\begin{equation}
\left\{
\begin{aligned}
	&\partial_t \u + \partial_x \v =0, \\
	&\partial_t \v + \frac{1}{\varepsilon^2} \partial_x \u = - \frac{1}{\varepsilon^2} \v.
\end{aligned}
\right. 
\label{I61}
\end{equation}
The original kinetic variables are recovered through relations 
\begin{equation}\label{eq:backtrasform}
f_+ = \dfrac{\u+\varepsilon \v}{2}, \qquad f_- = \dfrac{\u-\varepsilon \v}{2}.
\end{equation}
Throughout the manuscript, we will systematically switch back and forth between the systems~\eqref{I62} and~\eqref{I61}. The motivation is that~\eqref{I62} is more convenient for the development of the Monte Carlo solver whereas ~\eqref{I61} permits to compute the diffusion limit simply using the leading order term.

In fact, as $\varepsilon \to 0$ from the second equation in ~\eqref{I61} we obtain the local equilibrium
\[
\v= -\partial_x \u,
\]
which, using the first equation~\eqref{I61}, shows that the behavior of the solution, at least formally, is governed by the heat equation
\begin{equation}
 \partial_t \u= \partial_{xx}\u. 
 \label{I7b}
\end{equation}

\subsection{A standard Monte Carlo method for the GT model\label{sec:MC-GT}}


To construct the Monte Carlo scheme, we define an ensemble of $N$ particles $\left\{X_k(t),V_k(t)\right\}_{k=1}^N$, in which $X_k(t)$ represents the position and $V_k(t)$ the velocity of particle $k$ at time $t$.  For each $k$, the velocity $V_k(t)$ can only take two values, namely 
\[
\mathcal{V}^\epsi_{\pm}=\pm\frac{1}{\varepsilon}.
\] 
We will approximate the functions $f_+$ (respectively $f_-$), representing the density of particles with velocity $\mathcal{V}^\epsi_+$ (respectively $\mathcal{V}^\epsi_-$), by an empirical distribution
\begin{align}
  \mu_+(x,t) &= \frac{m_p}{N}\sum_{k=1}^{N} \, \delta(x-X_k(t))
\delta_{\textrm{Kr}}(\mathcal{V}^\epsi_+-V_k(t)), \label{particle_1}\\
  \mu_-(x,t) &= \frac{m_p}{N}\sum_{k=1}^{N} \, \delta(x-X_k(t))
\delta_{\textrm{Kr}}(\mathcal{V}^\epsi_--V_k(t)),\label{particle_2}
\end{align}
in which $\delta$ represents a Dirac delta, $\delta_{\textrm{Kr}}$ a Kronecker delta, and the mass $m_p$ of an individual particle is defined as 
\begin{equation}\label{mass}
m_p=\frac{1}{N}\int_{\Omega}u(x,0)dx.
\end{equation}
Note that one can introduce a (possibly time-dependent) weight $w_k(t)$ to each particle, and consider the corresponding weighted empirical distributions. We will not pursue this path here, as doing so does not affect the modified Monte Carlo schemes that are the focus of this manuscript.
\begin{remark}
[From empirical distributions to space-discretized particle densities]
Given the ensemble of particles $\left\{X_k(t),V_k(t)\right\}_{k=1}^N$, a  approximation to the particle densities $f_{\pm}(x,t)$ can be obtained as an histogram by introducing a spatial mesh with centers $\left\{x_j\right\}_{j=1}^J$ and mesh spacing $\Delta x$, and defining the stochastic approximation $f_{\pm}(x_j,t)$ by simply counting the particles inside the bin $[x_{j}-\Delta x/2,x_j+\Delta x/2]$:
\begin{equation}\label{eq:histogram}
f_{\pm}(x_j,t)=\int_{x_j-\Delta x/2}^{x_j+\Delta x/2}1\cdot d\mu^{\pm}(x,t).
\end{equation}
The density $\u(x_j,t)$ can then be obtained as 
\begin{equation}
\u(x_j,t)=f_{+}(x_j,t)+f_{-}(x_j,t).
\end{equation}
Clearly, many alternative approaches are available in the literature, see, e.g., \cite{LP, sliu}. In particular, a popular alternative to the above-described histogram approach is kernel density estimation \cite{scottKDE}.
\end{remark}

To simulate the Goldstein-Taylor model, one needs to define a stochastic process for the evolution of the ensemble $\left\{X_k(t),V_k(t)\right\}_{k=1}^N$, such that the population dynamics corresponds to~\eqref{I62}.  We introduce the discretized time $t^n=n\Delta t$, with $\Delta t$ the time step and $n\ge 0$, and write the time-discretized particle states as $X_k^n\approx X_k(t^n)$ and $V_k^n\approx V_k(t^n)$ respectively \cite{Seaid}. 
The standard Monte Carlo method for the Goldstein-Taylor model \eqref{I62} is based on a splitting between the transport and collision terms \cite{LP, Seaid}
\begin{enumerate}
	\item{\textbf{Transport}:
\begin{equation}\label{eq:2speed-transport}
\begin{aligned}
&\partial_t f_+ + \dfrac{1}{\varepsilon}\partial_x f_+ = 0,\\
&\partial_t f_- - \dfrac{1}{\varepsilon}\partial_x f_- = 0,
\end{aligned}
\end{equation}
}
\item{\textbf{Collision}:
\begin{equation}\label{eq:2speed-collision}
\begin{aligned} 
&\partial_t f_+= \dfrac{1}{\varepsilon^2}\left(\dfrac{\u}{2}- f_+\right),\\
&\partial_t f_-= \dfrac{1}{\varepsilon^2}\left(\dfrac{\u}{2}- f_-\right).
\end{aligned}
\end{equation}
}
\end{enumerate}
The splitting~\eqref{eq:2speed-transport}--\eqref{eq:2speed-collision} provides a convenient strategy to define the evolution of the particles: the transport step can be seen to affect the particle positions, leaving the velocities untouched, whereas the collision step updates the velocities, leaving the positions unaltered.

\paragraph{Transport step} During the transport step (\ref{eq:2speed-transport}), each particle  advances from time $t^n$ over a time step of size $\Delta t$ by changing its position according to 
\begin{equation}
X_k^{n+1}=X_k^n+V_k^n\Delta t.\label{transport}
\end{equation}
This can easily be shown by inserting (\ref{particle_1})--(\ref{particle_2}) in
(\ref{eq:2speed-transport}). Thus, after the transport step, we obtain the intermediate empirical distributions,
\begin{align}
  \tilde{\mu}_+^{n}(x) &= \frac{m_p}{N}\sum_{n=1}^{N} \, \delta(x-X_k^{n+1})
\delta_{\textrm{Kr}}(\mathcal{V}^\epsi_+-V_k^n), \label{intermediate_particle_1}\\
  \tilde{\mu}_-^{n}(x) &= \frac{m_p}{N}\sum_{n=1}^{N} \, \delta(x-X_k^{n+1})
\delta_{\textrm{Kr}}(\mathcal{V}^\epsi_--V_k^n),\label{intermediate_particle_2}
\end{align}
from which the intermediate particle densities $\tilde{f}_{\pm}^{n}(x_j)$ can be computed via \eqref{eq:histogram}.

\paragraph{Collision step} We consider now the solution of the collision step (\ref{eq:2speed-collision}). First observe that, during the collision step, the density $\u$ is constant. Hence, the density after the full time step has already been obtained during the transport step, and we have $\u^{n+1}=\tilde{\u}^n$. The effect of the collision step is thus to randomly change the velocities of a fraction of the particles from $\mathcal{V}^\epsi_+$ to $\mathcal{V}^\epsi_-$, and vice versa, without affecting their positions. Thus, the density after the complete time step is given as 
\begin{equation}
\u^{n+1}(x)= \tilde{\u}^{n}(x)=\tilde{f}_{+}^{n}(x)+\tilde{f}_-^{n}(x),
\end{equation}
and its stochastic approximation follows from equation~\eqref{eq:histogram}.

Now observe that, since equation~(\ref{eq:2speed-collision}) is linear and local, its exact solution is known. When neglecting (for now) the stochastic particle discretizations~\eqref{eq:histogram}, and using as initial conditions the values of $f_+$ and $f_-$ after the transport step, i.e. $\tilde{f}_{\pm}^{n}(x)$, we have
 \begin{equation}\label{eq:2speed-collision_sol}
 \begin{aligned} 
 &f_+^{n+1}(x)= \exp(-\Delta t/\varepsilon^2)\tilde{f}_{+}^{n}(x)+\left(1-\exp(-\Delta t/\varepsilon^2)\right) \dfrac{\u^{n+1}(x)}{2},\\
 &f_-^{n+1}(x)= \exp(-\Delta t/\varepsilon^2)\tilde{f}_{-}^{n}(x)+\left(1-\exp(-\Delta t/\varepsilon^2)\right) \dfrac{\u^{n+1}(x)}{2}.
 \end{aligned}
 \end{equation}
At the Monte Carlo level, the above formula can be interpreted as the convex combination of two probability distributions: 
\begin{itemize}
  \item With probability $\exp(-\Delta t/\varepsilon^2)$, the speed of a particle does not change, and the particle is left untouched.
  \item With probability $\left(1-\exp(-\Delta t/\varepsilon^2)\right)$ the speed of a particle changes, and the new velocity is chosen to be 
 $\mathcal{V}^\epsi_+$ or $\mathcal{V}^\epsi_-$ with equal probability.
\end{itemize}




\begin{remark}[Computational complexity in the diffusion limit]\label{rem:Eff_GT}
Even if there is no time step restriction in the collision part of the algorithm (the step always represents a convex combination of two distributions, regardless of the size of $\Delta t$ and $\varepsilon$), the approach outlined above suffers of severe time step restrictions when $\varepsilon\rightarrow 0$.
In fact, the main problem arises in the transport phase of the algorithm: when $\varepsilon \to 0$, the scaled particle velocities diverge, since ${\mathcal V}^\epsi_+ \to \infty$ and ${\mathcal V}^\epsi_- \to -\infty$. This means that we are forced to take $\Delta t=O(\varepsilon)$, making the Monte Carlo solver unusable for small values of $\varepsilon$ in the diffusion limit. 
\end{remark}

While the standard Monte Carlo method suffers from a prohibitive computational cost, we can observe that in the limit, the kinetic equation becomes equivalent to the heat equation (\ref{I7b}) and consequently a Monte Carlo method which solves this problem is easily available. In this case, the Monte Carlo method consists in first assigning the positions to $N$ particles
which approximates the function $u(x,t=0)=u^0(x)$ at the initial time by the empirical measure $\mu_u^0(x)$
\begin{equation}
\mu_u^0(x) = \frac{m_p}{N}\sum_{n=1}^{N} \, \delta(x-X_k^0), \label{particle_diffusion}\\
\end{equation}
where the particle positions $X_k^0$ are sampled from the probability distribution with density $u^0(x)/m_p$ and the constant $m_p$ is defined as
\begin{equation}
m_p=\frac{1}{N}\int_{\Omega}u^0(x)dx,
\end{equation}
Successively, the position of the particles evolves in time according to
\begin{equation}\label{scheme:diff}
 X_k(t^n+\Delta t)=X_k(t^n)+\sqrt{2\Delta t}\xi_k^n,
\end{equation}
where $\xi_k^n \sim \mathcal{N}(0,1)$ is a standard normally distributed random number. Observing that the Monte Carlo method~\eqref{scheme:diff} does not have time step limitations, one would like to construct a Monte Carlo scheme for the Goldstein-Taylor model~\eqref{I61} which in the limit $\varepsilon\rightarrow 0$ degenerates to (\ref{scheme:diff}) without time step limitations induced by the characteristic speeds. This property is refereed to as Asymptotic-Preserving (AP) property and we discuss such a Monte Carlo scheme in the next paragraph.

\subsection{An AP Monte Carlo method for the GT model}\label{AP-GT}
In this section, we introduce a new Monte Carlo approach based on a suitable reformulation of the original system.
We first discuss a time discretization that leads to a reformulated GT model (section~\ref{sec:ref-GT}). Subsequently, we introduce our new scheme based on the reformulated system (section~\ref{sec:ap-mc-gt}).

\subsubsection{An implicit time-discrete reformulation\label{sec:ref-GT}}
We start considering the following fully implicit discretization for system (\ref{I61}) 
\begin{equation}
 \label{eq:SP1}
\begin{cases}\begin{aligned}
&\dfrac{\u^{n+1}-\u^n}{\dt} = - \partial_x \v^{n+1}, \\
&\varepsilon^2\dfrac{\v^{n+1}-\v^n}{\dt} = -\left(\partial_x \u^{n+1} + \v^{n+1}\right). 
\end{aligned}
\end{cases}
\end{equation}
Now, solving the second equation for $\v^{n+1}$, one obtains 
\begin{equation}\label{eq:help1}
    \v^{n+1} = \frac{\epsi^2}{\epsi^2+\dt}\v^n - \frac{\dt}{\epsi^2+\dt}\partial_x\u^{n+1},
\end{equation}
or, equivalently,
\begin{equation}
\label{eq:second}
    \dfrac{\v^{n+1}-\v^n}{\Delta t} + \frac{1}{\epsi^2+\dt}\partial_x\u^{n+1} = -\frac{1}{\epsi^2+\dt}\v^n.
\end{equation}
Plugging~\eqref{eq:help1} in the first equation we get 
\[
\begin{aligned}
\frac{\u^{n+1}-\u^n}{\dt} + \frac{\epsi^2}{\epsi^2+\dt}\partial_x\v^n & = \frac{\dt}{\epsi^2+\dt}\partial_{xx}\u^{n+1}.
\end{aligned}
\]
Finally, using the first equation of \eqref{eq:SP1} and filling this into (\ref{eq:second}) we get the equivalent form
\begin{equation}
 \label{eq:SP2n_GS}
\begin{cases}
\begin{aligned}
&\frac{\u^{n+1}-\u^n}{\dt} + \frac{\epsi^2}{\epsi^2+\dt}\partial_x\v^n = \frac{\dt}{\epsi^2+\dt} \partial_{xx}\u^{n+1},
\\
&\frac{\v^{n+1}-\v^n}{\dt} + \frac{1}{\epsi^2+\Delta t}\partial_{x}\u^n  = - \frac{1}{\epsi^2+\Delta t}\v^{n}+\dfrac{\Delta t}{\epsi^2+\dt}\partial_{xx}\v^{n+1},
\end{aligned}
\end{cases}
\end{equation}
which, using the change of variables~\eqref{eq:backtrasform}, can be also written as
\begin{equation}\label{eq:bgk-modified-2speed-disc}
\begin{cases}
\begin{aligned}
&\frac{f_{+}^{n+1}-f_{+}^{n}}{\Delta t}  + \dfrac{\varepsilon}{\varepsilon^2+\Delta t}\partial_x f_+^{n} =\dfrac{\Delta t}{\varepsilon^2+\Delta t}\partial_{xx}f_+^{n+1} + \dfrac{1}{\varepsilon^2+\Delta t}\left(\dfrac{\u^n}{2}-f_+^{n}\right)\\
&\frac{f_{-}^{n+1}-f_{-}^{n}}{\Delta t}  - \dfrac{\varepsilon}{\varepsilon^2+\Delta t}\partial_x f_-^{n} = \dfrac{\Delta t}{\varepsilon^2+\Delta t}\partial_{xx}f_-^{n+1}+ \dfrac{1}{\varepsilon^2+\Delta t}\left(\dfrac{\u^n}{2}-f_-^{n}\right).
\end{aligned}
\end{cases}
\end{equation}
Observe now that, by using \[
\displaystyle \frac{\u^{n+1} - \u^n}{\Delta t} = \partial_t \u + \mathcal{O}(\Delta t), \quad \frac{\v^{n+1} - \v^n}{\Delta t} = \partial_t \v + \mathcal{O}(\Delta t),
\]
we obtain that system (\ref{eq:SP2n_GS}) is equivalent up to first order in $\dt$ to
\begin{equation}
 \label{eq:SP2bis}
 \begin{cases}
\begin{aligned}
&\partial_t\u + \frac{\epsi^2}{\epsi^2+\dt}\partial_x\v  = \frac{\dt}{\epsi^2+\dt}\partial_{xx}\u ,\\
&\partial_t\v + \frac{1}{\epsi^2+\Delta t}\partial_x \u  = - \frac{1}{\epsi^2+\Delta t}\v + \frac{\dt}{\epsi^2+\dt}\partial_{xx}\v,
\end{aligned}
\end{cases}
\end{equation}
and in diagonal form  
\begin{equation}\label{eq:bgk-modified-2speed}
\begin{cases}
\begin{aligned}
&\partial_t f_+ + \dfrac{\varepsilon}{\varepsilon^2+\Delta t}\partial_x f_+ =\dfrac{\Delta t}{\varepsilon^2+\Delta t}\partial_{xx}f_+ + \dfrac{1}{\varepsilon^2+\Delta t}\left(\dfrac{\u}{2}-f_+\right)\\
&\partial_t f_- - \dfrac{\varepsilon}{\varepsilon^2+\Delta t}\partial_x f_- = \dfrac{\Delta t}{\varepsilon^2+\Delta t}\partial_{xx}f_- + \dfrac{1}{\varepsilon^2+\Delta t}\left(\dfrac{\u}{2}-f_-\right).
\end{aligned}
\end{cases}
\end{equation}
Note that, the left part of system (\ref{eq:bgk-modified-2speed}) is hyperbolic with characteristic speeds
\begin{equation}
\label{char}
\lambda^{\pm}(\Delta t,\varepsilon) = \pm\dfrac{\varepsilon}{\varepsilon^2+\Delta t}.
\end{equation}
When $\dt\to 0$ for a fixed $\epsi$, system (\ref{eq:SP2bis})
converges to the original system (\ref{I61}), while the characteristic speeds converge to the usual ones, i.e., \[\lambda^{\pm}(0,\varepsilon) =\pm \frac{1}{\varepsilon}.\]
On the other hand, for a fixed $\dt$, the characteristic speeds $\lambda^+$ and $\lambda^-$ are bounded for any value of $\varepsilon$ and converge to zero as $\epsi\to 0$, while the diffusion coefficient tends to $1$ in that limit. Consequently, the system becomes fully parabolic and converges to the solution of the heat equation 
\begin{equation}
\partial_t \u=\partial_{xx}\u.
\end{equation}

\subsubsection{The APMC method\label{sec:ap-mc-gt}}

We now introduce a Monte Carlo scheme that solves the Goldstein-Taylor model (\ref{I61}) for all choices of the time step $\Delta t$ and $\varepsilon$, without any  $\varepsilon$-dependent time step restriction. The method is based on the following splitting approach between:

\begin{enumerate}
\item{\textbf{Transport--diffusion}:
\begin{equation}\label{eq:2speed-modified-transport-f}
\begin{cases} 
\begin{aligned}
&\partial_t f_+ + \dfrac{\varepsilon}{\varepsilon^2+\Delta t}\partial_x f_{+} = \dfrac{\Delta t}{\varepsilon^2+\Delta t}\partial_{xx}f_{+},\\
&\partial_t f_- - \dfrac{\varepsilon}{\varepsilon^2+\Delta t}\partial_x f_{-} = \dfrac{\Delta t}{\varepsilon^2+\Delta t}\partial_{xx}f_{-}.
\end{aligned}
\end{cases}
\end{equation}
}
\item{\textbf{Collision}:
\begin{equation}\label{eq:2speed-modified-collision-f}
\begin{cases} 
\begin{aligned}
&\partial_t f_+ = \dfrac{1}{\varepsilon^2+\Delta t}\left(\dfrac{u}{2}- f_{+}\right), \\
&\partial_t f_- = \dfrac{1}{\varepsilon^2+\Delta t}\left(\dfrac{u}{2}- f_{-}\right). 
\end{aligned}
\end{cases}
\end{equation}
}
\end{enumerate}
Now, as we did with the standard Monte Carlo method in Section~\ref{sec:MC-GT}, we approximate the functions $f_+$ and $f_-$ by
a finite set of particles $\left\{X_k(t),V_k(t)\right\}_{k=1}^N$, which correspond to the empirical measures~\eqref{particle_1} and~\eqref{particle_2}. The velocities $V_k(t)$ now can take the two values
\[
{\mathcal V}^\epsi_{\pm}=\pm\frac{\varepsilon}{\varepsilon^2+\Delta t},
\]
which are bounded for any value of $\epsi$ and are such that ${\mathcal V}^\epsi_{\pm}\to 0$ as $\epsi\to 0$. 
Recall that we can then reconstruct a histogram approximation of the distributions $f^\pm$ on a mesh via~\eqref{eq:histogram}.

\paragraph{Transport--diffusion step} To solve the transport--diffusion step (\ref{eq:2speed-modified-transport-f}), we observe that the particle velocities now scale with ${\varepsilon}/({\varepsilon^2+\Delta t})$, and that the diffusion corresponds to a Brownian motion with coefficient $\Delta t/(\varepsilon^2+\Delta t)$.  Thus, particles move according to
\begin{equation}
X_k^{n+1}=X^{n}_k+{\Delta t} V^n_k+\sqrt{2\dfrac{\Delta t^2}{\varepsilon^2+\Delta t}}\xi_k^n,\qquad 1\le k \le N,\label{transport_modified}\end{equation}
in which $\xi_i^n \sim \mathcal{N}(0,1)$ is a standard normally distributed random variable. The velocity of the particles does not change in this step, and we again have the intermediate empirical distributions~\eqref{intermediate_particle_1} and~\eqref{intermediate_particle_2}, from which the intermediate particle densities $\tilde{f}_\pm^{n}(x_j)$ can be computed via~\eqref{eq:histogram}. 

\paragraph{Collision step}
We consider now the collision step (\ref{eq:2speed-modified-collision-f}). As in the standard Monte Carlo method in Section~\ref{sec:MC-GT}, the density $\u$ is constant. The effect of this step is to randomly change the velocities of a fraction of the particles from $\mathcal{V}_{\pm}^\epsi$ to $\mathcal{V}_{\mp}^\epsi$, keeping the positions untouched.  By solving (\ref{eq:2speed-modified-collision-f}) with the forward Euler method we get

\begin{equation}\label{eq:2speed-modified-collision-f-discr}
\begin{cases} 
\begin{aligned}
&\dfrac{f_+^{n+1}-\tilde{f}_+^{n}}{\Delta t} = \dfrac{1}{\varepsilon^2+\Delta t}\left(\dfrac{\u^{n+1}}{2}- \tilde{f}_+^{n}\right), \\
&\dfrac{f_-^{n+1}-\tilde{f}_-^{n}}{\Delta t} = \dfrac{1}{\varepsilon^2+\Delta t}\left(\dfrac{\u^{n+1}}{2}- \tilde{f}_-^{n}\right).
\end{aligned}
\end{cases}
\end{equation}
This forward Euler discretization leads to the following convex combination 
\begin{equation}\label{eq:2speed-modified-collision-f-solv-discr}
\begin{cases} 
\begin{aligned}
&f_+^{n+1}=\dfrac{\varepsilon^2}{\varepsilon^2+\Delta t}\tilde{f}_+^{n} +\frac{\Delta t}{\varepsilon^2+\Delta t}\dfrac{\u^{n+1}}{2}, \\
&f_-^{n+1}=\dfrac{\varepsilon^2}{\varepsilon^2+\Delta t}\tilde{f}_-^{n} +\frac{\Delta t}{\varepsilon^2+\Delta t}\dfrac{\u^{n+1}}{2}. \\
\end{aligned}
\end{cases}
\end{equation}
Compared with collision step~\eqref{eq:2speed-collision_sol} in the standard Monte Carlo method, the only change is a slightly changed collision rate, i.e., a slightly different probability of a velocity change.
At the Monte Carlo level, equation (\ref{eq:2speed-modified-collision-f-solv-discr}) is thus again interpreted as 
\begin{itemize}
  \item With probability ${\varepsilon^2}/({\varepsilon^2+\Delta t})$, the speed of a particle does not change, and the particle is left untouched.
  \item With probability ${\Delta t}/({\varepsilon^2+\Delta t})$ the speed of a particle changes, and the new velocity is chosen to be 
 $\mathcal{V}^\epsi_+$ or $\mathcal{V}^\epsi_-$ with equal probability.
\end{itemize}

It is easy to see the asymptotic-preserving property of the new method. In fact, the time step of the transport--diffusion step is now independent of $\varepsilon$. In particular, in the limit $\varepsilon\to 0$ we get a standard Brownian motion for the heat equation
\begin{equation}
X_k^{n+1}=X^{n}_k+\sqrt{2\Delta t}\xi_k^n,\qquad 1\le k \le N.
\label{diffusion}
\end{equation}

\begin{remark}
\begin{itemize}
\item The Goldstein-Taylor model is equivalent to the {telegrapher's equation}, other probabilistic approach can be derived using this latter form \cite{AR,Kac74}. Note however, that in the diffusion limit the time step in the above approaches has to be taken of the size of $\epsi$ and therefore the methods are not asymptotic-preserving.   
\item For deterministic solvers, other form of asymptotic-preserving splitting for the diffusion limit have been proposed \cite{BPR12, JPT1, JPT2, Klar, LemMieu2008,NP2}. However, these splitting do not possess a clear probabilistic interpretation and it is not immediate to use them in a Monte Carlo setting. 
\end{itemize}
\end{remark}


\section{Asymptotic-preserving Monte Carlo methods for the radiative transport}\label{R-T}

In this section we show how to generalize the above approach to the radiative transport equation \cite{Case, Cha} under the diffusive scaling. 

\subsection{The radiative transport equation}
Let $f(x,v,t)$ be the probability density distribution for particles at space point $x \in \RR^{d_x}$, at time $t$ traveling in direction $v \in \Omega \subseteq \RR^{d_v}$, with $\int_\Omega dv = S$. Particles undergo two types of interactions: scattering, with scattering coefficient $\sigma_s(x)$ and absorption, with absorption coefficient $\sigma_a(x)$. Under the diffusive scaling, $f$ solves the radiative transfer equation
\begin{equation}\label{radiative_t}
\partial_t f + \dfrac{v}{\varepsilon} \cdot\nabla_x f =\frac{1}{\varepsilon^2}\left(\sigma_s\rho- \sigma f\right)+ G,
\end{equation}
where $\sigma(x) = \sigma_s (x) +\epsi^2 \sigma_a(x)$ is the total transport coefficient, $G(x)$ is the source term, $\varepsilon>0$ is proportional to the mean free path and 
\begin{equation}\label{eq:RT-density}
\rho(x,t)=\frac{1}{S}\int_\Omega f dv'
\end{equation}
is the position density.


To study the process in the diffusive limit when $\varepsilon$ tends to zero, we use the following expansion in $\varepsilon$ of the distribution function
\begin{equation}\label{expansion}
f = f^{(0)} + \varepsilon f^{(1)} + \varepsilon^2 f^{(2)}+\ldots
\end{equation}
and we introduce it in (\ref{radiative_t}). Then, considering terms of the same order in $\varepsilon$, we get at the leading order
\begin{equation}
 f^{(0)}(x,v,t)=\frac{1}{S}\int_\Omega f^{(0)} dv'=\rho(x,t),
\end{equation}
where $\rho(x,t)$ represents the density of the gas. Then, to the first order in
$\varepsilon$, we get
\begin{equation}\label{second_order_eps}
 v\cdot \nabla_x f^{(0)}=\frac{\sigma_s}{S}\int_\Omega f^{(1)}dv'-\sigma_s f^{(1)}.
\end{equation}
Now, writing the balance equation in terms of $\varepsilon$ of  (\ref{radiative_t}), one gets
\begin{equation}
 \partial_t \rho+v\cdot \nabla_x f^{(1)}=\frac{\sigma_s}{S}\int_\Omega f^{(2)}dv'-\sigma_s f^{(2)}-\sigma_a \rho+G
\end{equation}
and the integration in velocity space yields 
\begin{equation}
 \partial_t \rho+\frac{1}{S}\int_\Omega (v\cdot \nabla_x f^{(1)})dv=-\sigma_a\rho+G.
\end{equation}
Now, assuming $\sigma_s$ positive and strictly bounded away from zero, since one has from (\ref{second_order_eps}) that
\begin{equation}
 f^{(1)}=-\frac{1}{\sigma_s} v\cdot \nabla_x \rho+\frac{1}{S}\int_\Omega f^{(1)}dv'= -\frac{1}{\sigma_s} v\cdot \nabla_x \rho+\rho^{(1)},
\end{equation}
we finally obtain the following equation 
\begin{equation}\label{heat_eq}
 \partial_t \rho=\frac{1}{S}\int_\Omega v\cdot \nabla_x\left(\frac{v}{\sigma_s}\cdot\nabla_x \rho\right)dv-\sigma_a\rho+G=D\nabla_x\cdot\left(\frac{1}{\sigma_s}\nabla_x \rho\right)-\sigma_a\rho +G,
\end{equation}
where $D$ is the so-called diffusion coefficient which, for example, takes the value $D=1/3$ in one-dimensional slab geometry and $D=1/2$ when $\Omega$ is a unit circle in two dimensions. 

\subsection{A standard Monte Carlo scheme for the radiative transport\label{sec:standard-MC-RT}}

Let us now discuss a standard Monte Carlo method for solving the radiative transport equation, highlighting the limitations
of this approach when close to the diffusive limit. The starting point is, as for the two-speed case, a time splitting scheme for (\ref{radiative_t}) (where for simplicity we set the source term $G(x)=0$). It reads 
\begin{enumerate}
\item{\textbf{Transport}: 
\begin{equation}\label{eq:transport_radiat}
\partial_t f + \dfrac{1}{\varepsilon}v\cdot \nabla_x = 0
\end{equation}
}
\item{\textbf{Collision}: 
\begin{equation}\label{eq:collision_radiat}
\partial_t f= \dfrac{\sigma_s}{\varepsilon^2}\left(\rho- f\right).
\end{equation}
}
\item{\textbf{Absorption}: 
\begin{equation}\label{eq:absorption_radiat}
\partial_t f= -\sigma_a f.
\end{equation}
}
\end{enumerate}

We again approximate the distribution by an empirical distribution,  using a finite set of particles with positions and velocities $\left\{X_k(t),V_k(t)\right\}_{n=1}^N$. The particle velocities are given as  
\[
V_k(t)=\frac{\tilde{V}_k(t)}{\varepsilon},\qquad \tilde{V}_k\in \Omega.
\]
Defining the mass $m_p$ of an individual particle as
\begin{equation}\label{mass_radiat}
m_p=\frac{1}{N}\int_{\mathbb{R}}\int_{\Omega}f(x,v,t=0)dvdx,
\end{equation}
we obtain the empirical distribution 
\begin{equation}\label{particle_radiat}
\mu(x,v,t) = \frac{m_p}{N}\sum_{n=1}^{N} \, \delta(x-X_k(t))\delta(v-V_k(t)),
\end{equation}
in which $\delta$ is again the Dirac delta.

\begin{remark}[Empirical particle densities in phase space]
We restrict for simplicity to the one-dimensional case in both $x$ and $v$. In analogy with the two-speed case, we can introduce a mesh and compute a histogram on the phase space mesh with cell centers $x_j$ and $v_{\ell}$ and mesh widths $\Delta x$ and $\Delta v$ as 
\begin{equation}\label{eq:hist_RT}
f(x_j,v_{\ell},t)=\int_{x_j-\Delta x/2}^{x_j+\Delta x/2}
\int_{v_{\ell}-\Delta v/2}^{v_{\ell}+\Delta v/2}
1\cdot d\mu(x,v,t).
\end{equation}
An empirical position density is then obtained as 
\begin{equation}\label{eq:hist_rho_RT}
\rho(x_j,t)=\int_{x_j-\Delta x/2}^{x_j+\Delta x/2}
\int_\Omega
1\cdot d\mu(x,v,t).
\end{equation}
\end{remark}

We can now describe the Monte Carlo method that correspond to transport and collisions of the particles. In the sequel we will neglect the presence of the source term $G$ which can be easily included in the method.

\paragraph{Transport}
As in the two-speed case, each particle advances from time $t^n$ over a time interval of length $\Delta t$ during the transport step (\ref{eq:transport_radiat}) by changing its position according to
\begin{equation}\label{transport_r}
X_k^{n+1}=X_k^n+V_k^n\Delta t.
\end{equation}
This can be shown by inserting (\ref{particle_radiat}) inside
(\ref{eq:transport_radiat}). 

We then have an intermediate empirical distribution:
\begin{equation}\label{particle_radiat_intermediate}
\tilde{\mu}^n(x,v) = \frac{m_p}{N}\sum_{n=1}^{N} \, \delta(x-X_k^{n+1})\delta(v-V_k^n),
\end{equation}
from which the intermediate particle density $\tilde{f}^n(x_j,v_\ell)$ can be computed using~\eqref{eq:hist_RT}.

\paragraph{Collision}
Next, we solve the collision process (\ref{eq:collision_radiat}) without absorption
\begin{equation}\label{eq:collision_radiat1}
\partial_t f= \dfrac{\sigma_s}{\varepsilon^2}\left(\rho- f\right),
\end{equation}
which corresponds to 
\begin{eqnarray}
{\tilde f}^{n+1}(x,v)&=& \exp\left(-\frac{\sigma_s\Delta t}{\varepsilon^2}\right)\tilde{f}^{n}(x,v)+\left(1-\exp\left(-\frac{\sigma_s\Delta t}{\varepsilon^2}\right)\right) \tilde{\rho}^{n}(x)\\
 \end{eqnarray}
 At the Monte Carlo level, the above formula can be interpreted in the following
way: 
\begin{itemize}
  \item With probability $\exp\left(-\sigma_s\Delta t/\varepsilon^2\right)$, the speed of a particle does not change
  \item With probability $\left(1-\exp\left(-\sigma_s\Delta t/\varepsilon^2\right)\right)$, the speed of a particle changes to a new value $V_k=\tilde{V}_k/\varepsilon$, in which $\tilde{V}_k$ is a random value with uniform probability in the domain $\Omega$.
\end{itemize}
  
\paragraph{Absorption}
We consider now the solution of the absorption step (\ref{eq:absorption_radiat}). Unlike the two-speed case, due to absorption the density $\rho(x,t)$ is not conserved. 

In a time step $\Delta t$, we solve the absorption process
\begin{equation}\label{eq:collision_abs}
\partial_t f= -{\sigma_a}f,
\end{equation}
which allows to compute  
\begin{eqnarray}
{f}^{n+1}(x,v)&=& \exp\left(-{\sigma_a\Delta t}\right)\tilde{f}^{n+1}(x,v),
 \end{eqnarray}
The above process is easily realized, assuming that with probability $1-\exp\left(-\sigma_a\Delta t\right)$, the particle gets absorbed and disappears from the simulation.
  
\begin{remark}[Mesh-based approach]\label{rem:mesh}
We may consider a method based on a mesh in space. We define the density of the particles in the center of the cells, $\rho(x_j,t)=(\int_\Omega f(x_j,v,t) dv')/S$ and solve (\ref{eq:collision_radiat}) in $x_j, \ j=1,..,M$ with $M$ the number of mesh points. 
In order to compute the integral of the distribution function in the cell centers different techniques can be used. The simplest first order space reconstruction in one dimension, the same used for the two speed case, is given by~\eqref{eq:hist_rho_RT}.
\end{remark}

Concerning the asymptotic behavior for small values of $\epsi$ the same remark can be made as in the two-speed case (see Remark~\ref{rem:Eff_GT}): in the limit when $\varepsilon$ tends to zero, the main computational bottleneck is due to transport phase, where the transport speeds approach infinity and hence infinitely small time steps would be required.  

Also for the radiative transfer equation, a standard Monte Carlo method for the diffusion equation~\eqref{heat_eq} can be derived. In the one dimensional case, if we neglect the source term and the absorption, it consists in initializing the system by creating an ensemble of particles $\left\{X_k^0\right\}_{n=1}^N$ that are sampled according to the local density $\rho(x,t=0)$,
and then by advancing in a time step $\Delta t$ following the equation
\begin{equation}\label{transport_heat}
X_k^{n+1}=X_k^n+\sqrt{2\dfrac{D}{\sigma_s}\Delta t}\xi_k^n,\end{equation}
where $\xi_k^n \sim \mathcal{N}(0,1)$ is a standard normally distributed random number. Thus, we want to construct am asymptotic-preserving  Monte Carlo method for radiative transfer~\eqref{radiative_t} that automatically degenerates to 
the above described Monte Carlo method without any time step restriction induced by the unbounded increasing particle speed when $\varepsilon\rightarrow 0$. 

\subsection{An AP Monte Carlo method for the radiative transport}\label{AP-RT}
In this section, we generalize the reformulation discussed in Section \ref{AP-GT} for the Goldstein-Taylor model~\eqref{I61} to the case of the radiative transfer model~\eqref{radiative_t}. 
We start by reformulating the radiative equation by using the even and odd formalism and by introducing a suitable time discretization (section~\ref{sec:ref-R-T}). In section~\ref{sec:ref-AP-RT}, we then make use of this reformulation for constructing our new scheme.

\subsubsection{The reformulated radiative transport equation\label{sec:ref-R-T}} 
In order to emphasize the analogies with the Goldstein--Taylor model we consider radiative transport equation~\eqref{radiative_t} without 
source term $G(x)=0$, i.e.,
\begin{equation}\label{radiative_t_1d}
 \partial_t f+\dfrac{v}{\varepsilon}\cdot \nabla_x f=\frac{1}{\varepsilon^2}\left(\sigma_s\rho-\sigma f\right)=\frac{\sigma_s}{\varepsilon^2}\left(\rho-f\right)-\sigma_a f
\end{equation}
with \[
\rho(x,t)=\dfrac{1}{S}\int_{\Omega} f(x,v',t) dv'.
\]

We first rewrite the radiative transfer equation as
\begin{equation}
\begin{cases}\begin{aligned}\label{radiative_t_split}
&\partial_t f_{+}+\dfrac{v}{\varepsilon}\cdot\nabla_x f_{+} =\frac{\sigma_s}{\varepsilon^2}\left(\rho-f_+\right)-\sigma_a f_+\\
 &\partial_t f_{-}-\dfrac{v}{\varepsilon}\cdot\nabla_x f_{-} =\frac{\sigma_s}{\varepsilon^2}\left(\rho-f_+\right)-\sigma_a f_-,
\end{aligned}
\end{cases}
\end{equation}
where now  $f_{+}(x,v,t)=f(x,v,t)$ and $f_{-}(x,v,t)=f(x,-v,t)$. 
This permits to define the even and odd parities 
\begin{equation}\label{eq:trasform-RT}
 r(x,v,t)=\frac{1}{2}\left(f_+(x,v,t)+f_-(x,v,t)\right), \qquad  j(x,v,t)=\frac{1}{2\varepsilon}\left(f_+(x,v,t)-f_-(x,v,t)\right).
\end{equation}
Then, $r(x,v,t)$ and $j(x,v,t)$ satisfy the following equivalent system
\begin{equation}
\begin{cases}\begin{aligned}
\label{even-odd}
 &\partial_t r+v\cdot\nabla_x j =\frac{\sigma_s}{\varepsilon^2}\left(\rho-r\right)-\sigma_a r\\
 &\partial_t j+\frac{v}{\varepsilon^2}\cdot\nabla_x r =-\frac{\sigma_s}{\varepsilon^2}j+\sigma_a j.
\end{aligned}
\end{cases}
\end{equation}
The inverse transformation of~\eqref{eq:trasform-RT} is easily seen to be
\begin{equation}\label{eq:backtrasform-RT}
f_+(x,v,t)=r(x,v,t)+\varepsilon j(x,v,t), \qquad f_-(x,v,t)=r(x,v,t)-\varepsilon j(x,v,t).
\end{equation}
In order to construct an implicit reformulation of the problem we first split the system into three parts as
\begin{equation}
{\rm (I)}
\begin{cases}\begin{aligned}
\label{even-odd-split1}
&\partial_t r+v\cdot\nabla_x j =0\\
&\partial_t j+\frac{v}{\varepsilon^2}\cdot\nabla_x r =-\frac{\sigma_s}{\varepsilon^2} j,
\end{aligned}
\end{cases}
\end{equation}
\begin{equation}
{\rm (II)}
\begin{cases}\begin{aligned}
\label{even-odd-split2}
&\partial_t r=\frac{\sigma_s}{\varepsilon^2}\left(\rho-r\right)\\
&\partial_t j=0
\end{aligned}
\end{cases}
\end{equation}
and
\begin{equation}
{\rm (III)}
\begin{cases}\begin{aligned}
\label{even-odd-split3}
&\partial_t r=-\sigma_a r\\
&\partial_t j=0.
\end{aligned}
\end{cases}
\end{equation}
The first step (I) now has the same structure of the Goldstein--Taylor model and we can follow the approach developed in Section~\ref{sec:ref-GT}. We consider the implicit discretization of (\ref{even-odd-split1}) as
\begin{equation}
\begin{cases}\begin{aligned}
\label{even-odd-disc}
 &\frac{r_*^{n+1}-r^n}{\Delta t}+v\cdot\nabla_x j_*^{n+1} =0\\
 &\frac{j_*^{n+1}-j^n}{\Delta t}+\frac{v}{\varepsilon^2}\cdot\nabla_x r_*^{n+1} =-\frac{\sigma_s}{\varepsilon^2}j_*^{n+1},
\end{aligned}
\end{cases}
\end{equation}
where $r_*^{n+1}$ and $j_*^{n+1}$ denote the solutions of this first step.

Solving for $j_*^{n+1}$ one gets
\begin{equation}\label{j_kew}
 j_*^{n+1}=\frac{\varepsilon^2}{\varepsilon^2+\sigma_s\Delta t} j^n-\frac{\Delta t}{\varepsilon^2+\sigma_s \Delta t}v\cdot\nabla_x r_*^{n+1},
\end{equation}
or, equivalently,
\begin{equation}\label{j_kew_help2}
 \dfrac{j_*^{n+1}-j^n}{\Delta t}+\frac{1}{\varepsilon^2+\sigma_s \Delta t}v\cdot\nabla_x r_*^{n+1}=-\frac{\sigma_s}{\varepsilon^2+\sigma_s \Delta t} j^n.
\end{equation}
Equation~\eqref{j_kew} can be plugged in the first equation of (\ref{even-odd-disc}) to give
\begin{equation}\label{r_kew}
 r_*^{n+1}=r^n-\Delta t v\cdot\nabla_x \left(\frac{\varepsilon^2}{\varepsilon^2+\sigma_s \Delta t} j^n-\frac{\Delta t}{\varepsilon^2+\sigma_s \Delta t}v\cdot\nabla_x r_*^{n+1}\right).
\end{equation}
Now, using the first equation of (\ref{even-odd-disc}) into (\ref{j_kew_help2}) gives
\begin{equation}\label{r-j_kew1}
\begin{cases}\begin{aligned}
&\frac{r_*^{n+1}-r^n}{\Delta t}+\frac{\varepsilon^2}{\varepsilon^2+\sigma_s \Delta t}v\cdot\nabla_x j^n =\frac{\Delta t}{\varepsilon^2+\sigma_s \Delta t} v\cdot \nabla_x\left(v\cdot\nabla_{x}r_*^{n+1}\right)\\
 &\frac{j_*^{n+1}-j^n}{\Delta t}+\frac{1}{\varepsilon^2+\sigma_s \Delta t}v\cdot\nabla_x r^{n}=\frac{\Delta t}{\varepsilon^2+\sigma_s \Delta t} v\cdot \nabla_x\left(v\cdot\nabla_{x}j_*^{n+1}\right)-\frac{\sigma_s}{\varepsilon^2+\sigma_s \Delta t}j^n.
\end{aligned}
\end{cases}
\end{equation}
The second part of the splitting, equation~\eqref{even-odd-split2}, can be discretized similarly with an implicit method to give
\begin{equation}\label{r-j_kew2}
\begin{cases}\begin{aligned}
&\frac{r_{**}^{n+1}-r_*^{n+1}}{\Delta t} =\frac{\sigma_s}{\varepsilon^2+\sigma_s\Delta t}\left(\rho_{**}^{n+1}-r_*^{n+1}\right)
\\
 &\frac{j_{**}^{n+1}-j_*^{n+1}}{\Delta t}=0,
\end{aligned}
\end{cases}
\end{equation}
where $\rho_{**}^{n+1}=\rho_{*}^{n+1}$ since the density remains unchanged during this step.
We observe now that (\ref{r-j_kew1})-(\ref{r-j_kew2}) are, up to an error $\mathcal{O}(\Delta t)$, equivalent to a time splitting of the reformulated system
\begin{equation}
\begin{cases}\begin{aligned}
\label{even-odd_reformulated}
 &\partial_t r+\frac{\varepsilon^2}{\varepsilon^2+\sigma_s\Delta t}
 v\cdot\nabla_x j =\frac{\Delta t}{\varepsilon^2+\sigma_s\Delta t} v\cdot \nabla_x\left(v\cdot\nabla_{x}r\right)+\frac{\sigma_s}{\varepsilon^2+\sigma_s\Delta t}\left(\rho-r\right)\\
 &\partial_t j+\frac{1}{\varepsilon^2+\sigma_s\Delta t}v\cdot\nabla_x r =\frac{\Delta t}{\varepsilon^2+\sigma_s\Delta t} v\cdot \nabla_x\left(v\cdot\nabla_{x}j\right)-\frac{\sigma_s}{\varepsilon^2+\sigma_s \Delta t}j.
\end{aligned}
\end{cases}
\end{equation}
Using the back transformation~\eqref{eq:backtrasform-RT}, equation~\eqref{even-odd_reformulated} can be written also as
\begin{equation}
\begin{cases}\begin{aligned}\label{radiative_t_split_reform}
& \partial_t f_{+}+\dfrac{\varepsilon}{\varepsilon^2+\sigma_s \Delta t}v\cdot\nabla_x f_{+} =\frac{\Delta t}{\varepsilon^2+\sigma_s \Delta t} v\cdot \nabla_x\left(v\cdot\nabla_{x}f_+\right)+\frac{\sigma_s}{\varepsilon^2+\sigma_s\Delta t}\left(\rho-f_+\right)\\
 &\partial_t f_{-}-\dfrac{\varepsilon}{\varepsilon^2+\sigma_s \Delta t}v \cdot\nabla_xf_{-} =\frac{\Delta t}{\varepsilon^2+\sigma_s \Delta t} v\cdot \nabla_x\left(v\cdot\nabla_{x}f_-\right)+\frac{\sigma_s}{\varepsilon^2+\sigma_s\Delta t}\left(\rho-f_-\right).
\end{aligned}
\end{cases}
\end{equation}
Note that, for fixed values of $\varepsilon$, the above equations revert to the original system (\ref{even-odd}) or (\ref{radiative_t_1d}) in the limit when the time step $\Delta t$ tends to zero when the absorption coefficient $\sigma_a=0$. Let observe that up to an error of order $\mathcal{O}(\Delta t)$ system (\ref{even-odd_reformulated}) or (\ref{radiative_t_split_reform}) plus the third step of the splitting (\ref{even-odd-split3}) represent a first order in time approximation of the original radiative transfer equation.
On the other hand, for all $\Delta t>0$, system (\ref{even-odd_reformulated}) or (\ref{radiative_t_split_reform}) together with (\ref{even-odd-split3}) are a $\mathcal{O}(\Delta t)$ approximation with bounded eigenvalues for every choice of $\Delta t$ of the original system. In particular, for every finite time step, the system tends to the limiting diffusion equation~\eqref{heat_eq} in the limit when $\varepsilon$ tends to zero.

\begin{remark}[Micro-macro decomposition]\label{rem:micro-macro}
As an alternative to the odd-even splitting above, one could also consider a micro-macro splitting, see, e.g., \cite{Klar,LemMieu2008}. Let us illustrate the approach in one space dimension and with $\sigma_s=\sigma_a=1$ for simplicity. In that case, we write 
\begin{equation}\label{eq:micro-macro-expansion}
	f(x,v,t) = \rho(x,t) + \varepsilon g(x,v,t), 
\end{equation}
with $\rho$ defined as before, from which we naturally derive that $\int_{-1}^1 g(x,v,t) dv = 0$. 
Inserting this expansion in~\eqref{radiative_t_1d} and averaging over velocity space, we get the system
\begin{equation}\label{eq:micro-macro-equation}
\begin{cases}
	&\partial_t \rho + \partial_x \langle vg\rangle = 0,\\
	&\partial_t g + \dfrac{1}{\varepsilon^2}v\partial_x \rho + \dfrac{1}{\varepsilon}\partial_x\left(vg-\langle vg\rangle\right)=-\dfrac{1}{\varepsilon^2}g,
\end{cases}
\end{equation}
in which we introduced the notation $\langle \cdot \rangle = (1/2)\int_{-1}^1\cdot \; dv$ to denote the average over velocity space. 
It can easily be checked that~\eqref{eq:micro-macro-equation} is equivalent to the original kinetic equation~\eqref{radiative_t_1d}. Following a similar reasoning as above, one can obtain up to ${\mathcal O}(\Delta t)$ the modified equation
\begin{equation}\label{eq:multispeed-IMEX-reform}
\begin{cases}\begin{aligned}	
	&	\partial_t\rho+\dfrac{\varepsilon^2}{\varepsilon^2+\Delta t}\partial_x \langle vg\rangle = \dfrac{\Delta t}{\varepsilon^2+\Delta t}\langle v^2\rangle\partial_{xx}\rho,\\
	&\partial_t g + \dfrac{1}{\varepsilon^2+\Delta t}v\partial_x \rho+\dfrac{\varepsilon}{\varepsilon^2+\Delta t}\partial_x\left(vg-\langle vg\rangle\right)=-\dfrac{1}{\varepsilon^2+\Delta t}g,
	\end{aligned}
\end{cases}
\end{equation}
Equation~\eqref{eq:multispeed-IMEX-reform} satisfies the same desirable properties as equation~\eqref{radiative_t_split_reform} or~\eqref{even-odd_reformulated}: it converges to the original equation~\eqref{radiative_t_1d} for fixed $\varepsilon$ as $\Delta t$ tends to zero, and to the limiting heat equation~\eqref{heat_eq} as $\varepsilon$ tends to zero for fixed $\Delta t$. Thus, equation~\eqref{eq:multispeed-IMEX-reform} may also serve as the basis for an asymptotic-preserving particle scheme, see Remark~\ref{rem:ap-micro-macro}. 
\end{remark}

\subsubsection{The APMC method\label{sec:ref-AP-RT}} 

In this paragraph, we show how the reformulation~\eqref{even-odd_reformulated} permits to develop a Monte Carlo scheme that is not limited by the stiffness of the equation~\eqref{radiative_t} in the limit when $\varepsilon$ tends to $0$. 
The Monte Carlo method is based on the following splitting of the reformulated system~\eqref{even-odd_reformulated}:
\begin{enumerate}
\item{\textbf{Transport and diffusion}:
\begin{equation}\label{eq:radiative-modified-transport-f}
\begin{cases} 
\begin{aligned}
&\partial_t f_+ + \dfrac{\varepsilon}{\varepsilon^2+\sigma_s\Delta t}v\cdot\nabla_x f_+ = \frac{\Delta t}{\varepsilon^2+\sigma_s\Delta t} v\cdot \nabla_x\left(v\cdot\nabla_{x}f_+\right),\\
&\partial_t f_- -\dfrac{\varepsilon}{\varepsilon^2+\sigma_s\Delta t}v\cdot\nabla_x f_- = \frac{\Delta t}{\varepsilon^2+\sigma_s\Delta t} v\cdot \nabla_x\left(v\cdot\nabla_{x}f_-\right).
\end{aligned}
\end{cases}
\end{equation}
}
\item{\textbf{Collision}:
\begin{equation}
\begin{cases} \label{eq:radiative-modified-collision-f}
\begin{aligned}
&\partial_t f_+ = \frac{\sigma_s}{\varepsilon^2+\sigma_s\Delta t}(\rho-f_+), \\
&\partial_t f_- = \frac{\sigma_s}{\varepsilon^2+\sigma_s\Delta t}(\rho-f_-). \end{aligned}
\end{cases}
\end{equation}
}
\item{\textbf{Absorption}:
\begin{equation}
\begin{cases} \label{eq:radiative-modified-absorption-f}
\begin{aligned}
&\partial_t f_+ = -\sigma_a f_+, \\
&\partial_t f_- = -\sigma_a f_-. \end{aligned}
\end{cases}
\end{equation}
}
\end{enumerate}

We are now ready to introduce the Monte Carlo method. We again approximate the distribution by an empirical distribution,  using a finite set of particles with positions and velocities $\left\{X_k(t),V_k(t)\right\}_{n=1}^N$, see also equation~\eqref{particle_radiat}. The particle velocities are now given as  
\[
V_k(t)=\dfrac{\varepsilon}{\varepsilon^2+\sigma_s\Delta t}\tilde{V}_k(t),\qquad \tilde{V}_k\in \Omega,
\]
and the mass of an individual particle is given by~\eqref{mass_radiat}. 

\paragraph{Transport and diffusion} The transport and diffusion step (\ref{eq:radiative-modified-transport-f}) can be handled by observing that (\ref{eq:radiative-modified-transport-f}) represents a population of particles each one moving according to
\begin{equation}
X^{n+1}_k=X^{n}_{n}+\Delta t V^n_k+\sqrt{2\dfrac{\Delta t^2 \left(V^n_k\right)^2}{\varepsilon^2+\sigma_s\Delta t}}\xi_k^n,\qquad 1\le n \le N,\label{transport_modified_rad}
\end{equation}
where $X^{n+1}_k$ indicates the new position of the particle after the transport and $\xi_k^n \sim \mathcal{N}(0,1)$ are independent standard normally distributed random numbers. The velocities of the particles do not change in this step, and we have the intermediate empirical distribution~\eqref{particle_radiat_intermediate}, from which an intermediate particle density $\tilde{f}^n(x_j,v_\ell)$ can be computed using for example~\eqref{eq:hist_RT} in the one-dimensional case.

\paragraph{Collision} We consider now the solution of the collision step (\ref{eq:radiative-modified-collision-f}). 
Let observe that the collision step preserves the density $\rho$, since we consider a situation without absorption.  Thus, collision does not affect particle positions.
We can write the solution of~\eqref{eq:radiative-modified-collision-f} using the original implicit formulation as
\begin{equation}\label{eq:mc-rt-coll}
\begin{cases} 
\begin{aligned}
&\tilde f_+^{n+1}=\dfrac{\varepsilon^2}{\varepsilon^2+\sigma_s\Delta t}\tilde{f}_+^{n} +\frac{\sigma_s\Delta t}{\varepsilon^2+\sigma_s\Delta t}{\tilde \u^{n+1}}, \\
&\tilde f_-^{n+1}=\dfrac{\varepsilon^2}{\varepsilon^2+\sigma_s\Delta t}\tilde{f}_-^{n} +\frac{\sigma_s\Delta t}{\varepsilon^2+\sigma_s\Delta t}{\tilde \u^{n+1}}. \\
\end{aligned}
\end{cases}
\end{equation}
At the Monte Carlo level, the above formulas can be interpreted in the following
way: 
\begin{itemize}
  \item With probability ${\varepsilon^2}/({\varepsilon^2+\sigma_s\Delta t})$, the speed of a particle does not change.
  \item With probability ${\sigma_s\Delta t}/({\varepsilon^2+\sigma_s\Delta t})$, the speed of a particle changes to a new value $V_k=\dfrac{\varepsilon}{\varepsilon^2+\sigma_s\Delta t}\tilde{V}_k$, in which $\tilde{V}_k$ is a random value with uniform probability in the domain $\Omega$.
\end{itemize}
\paragraph{Absorption} We consider now the solution of the absorption step (\ref{eq:radiative-modified-absorption-f}). This step is analogous to the step already discussed
in Section 3.2. However, instead of using the exact solution of equation (\ref{eq:radiative-modified-absorption-f}) to construct the Monte Carlo method, for consistency with the previous steps we use a first order implicit time discretization. This reads
\begin{equation}\label{eq:mc-rt-absor}
\begin{cases} 
\begin{aligned}
&f_+^{n+1}=\frac{1}{1+\sigma_a\Delta t}\tilde f^{n+1}_+, \\
&f_-^{n+1}=\frac{1}{1+\sigma_a\Delta t}\tilde f^{n+1}_-. \\
\end{aligned}
\end{cases}
\end{equation}
Let observe that this step modifies the total mass of the system, i.e. $\rho^{n+1}=\tilde \rho^{n+1}/({1+\sigma_a\Delta t})$.
At the Monte Carlo level, the above formulas can be simply interpreted as: with probability ${\sigma_a\Delta t}/({1+\sigma_a\Delta t})$ a particle is removed from the domain. 



When $\varepsilon\to 0$ and in the case without absorption, as for the case of the Goldstein-Taylor model, the scheme automatically reduces to a standard Monte Carlo scheme for the diffusion equation,
\begin{equation}
X^{n+1}_k=X^{n}_{n}+\sqrt{2{\Delta t \left(V^n_k\right)^2}}\xi_k^n,\qquad 1\le n \le N,
\end{equation}
i.e., the scheme satisfies the asymptotic preserving property. In fact, in this limit, the scheme degenerate to the solution of the first step in which the transport speed is zero. The relaxation step clearly does not play any role in this limit. In the case with absorption, in the same limit, the scheme degenerates to a Monte Carlo method for the diffusion reaction equation.

\begin{remark}[Derivation based on the micro-macro decomposition]\label{rem:ap-micro-macro}
One can also derive a particle scheme starting from the micro-macro decomposition~\eqref{eq:micro-macro-expansion}, see Remark~\ref{rem:micro-macro}.  The modified equation~\eqref{eq:multispeed-IMEX-reform} then leads to an asymptotic-preserving Monte Carlo scheme that is very similar to the scheme derived above. The only difference is that the transport and diffusion step~\eqref{transport_modified_rad} then reads
\begin{equation}
X^{n+1}_k=X^{n}_{n}+\Delta t V^n_k+\sqrt{2\dfrac{\Delta t^2 D}{\varepsilon^2+\Delta t}}\xi_k^n,\qquad 1\le n \le N,\label{transport_modified_rad_micro-macro}
\end{equation}
with $D=\langle v^2\rangle$. In~\eqref{transport_modified_rad_micro-macro}, the diffusive correction is performed with the diffusion coefficient of the limiting heat equation~\eqref{heat_eq}, and not with a diffusion coefficient that depends on the velocity of the particle.  
\end{remark}

\section{Numerical experiments}\label{Num}
In this section we discuss several numerical tests with the aim of demonstrating the behavior and the performance of the new Asymptotic Preserving Monte Carlo scheme detailed in the previous
Sections. We start by discussing the case of the Goldstein-Taylor model and we end with the radiative transport problem. For all the tests considered we compare our scheme with a finite volume 
asymptotic preserving method, the one described in \cite{BPR17}. The reference method is based on an implicit-explicit time discretization which makes it unconditionally stable with respect the scaling parameter $\varepsilon$, the sole stability conditions being dictated by the diffusive or the hyperbolic regime. 

\subsection{The Goldstein-Taylor model}
\begin{figure}[ht!]
\begin{center}
    \includegraphics[scale=0.37]{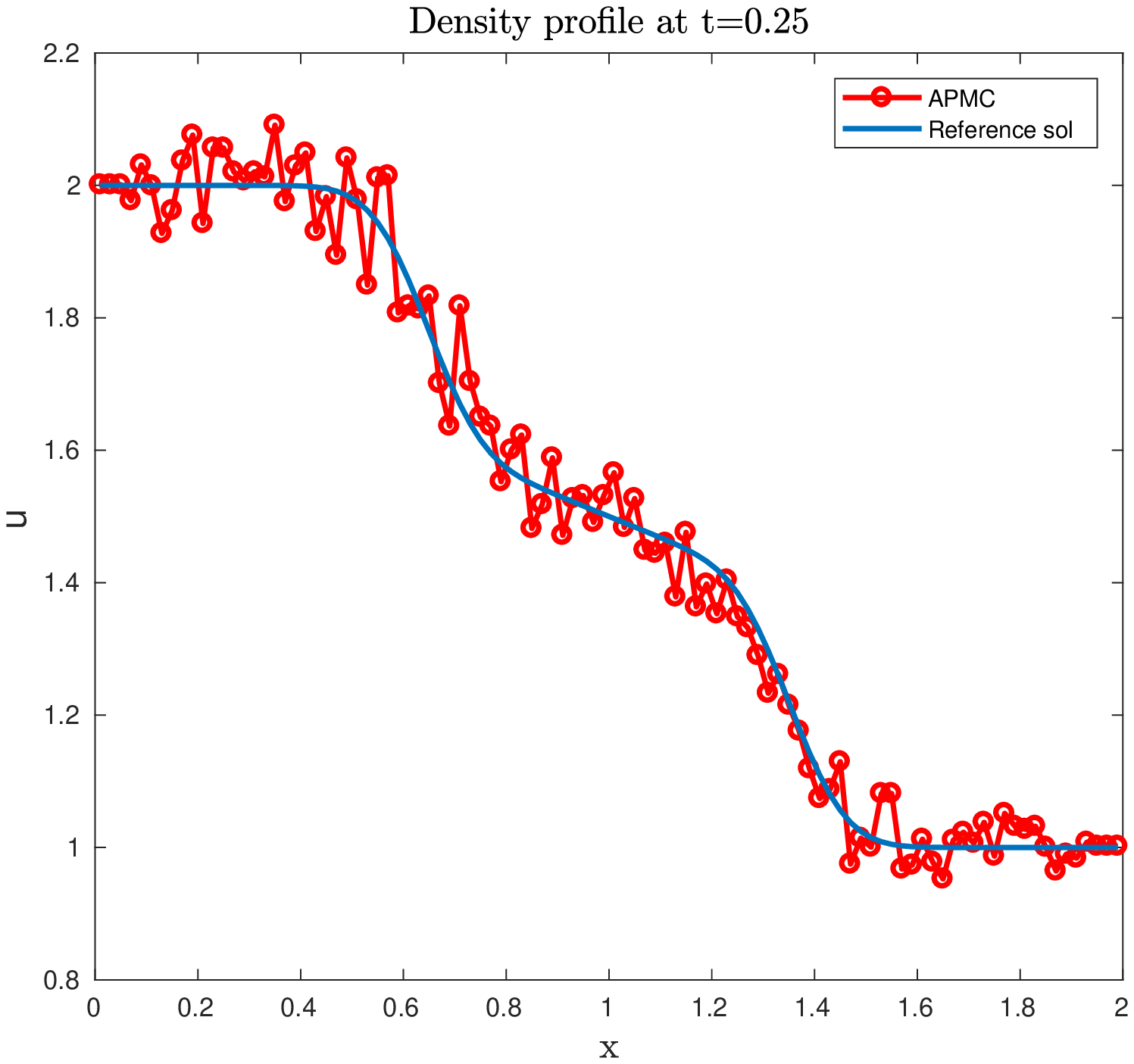}\hspace{1cm}
        \includegraphics[scale=0.37]{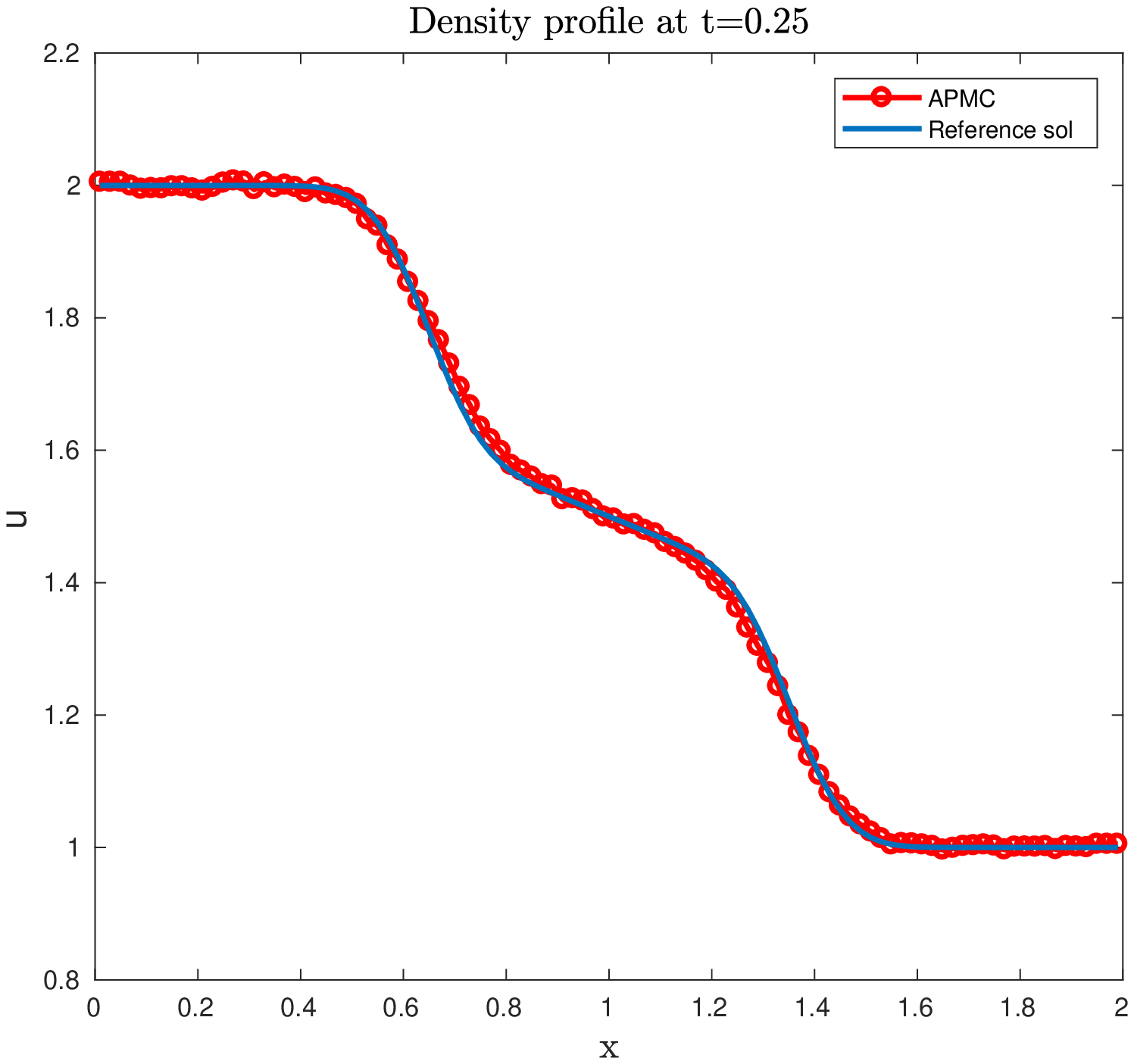}\\
     \includegraphics[scale=0.37]{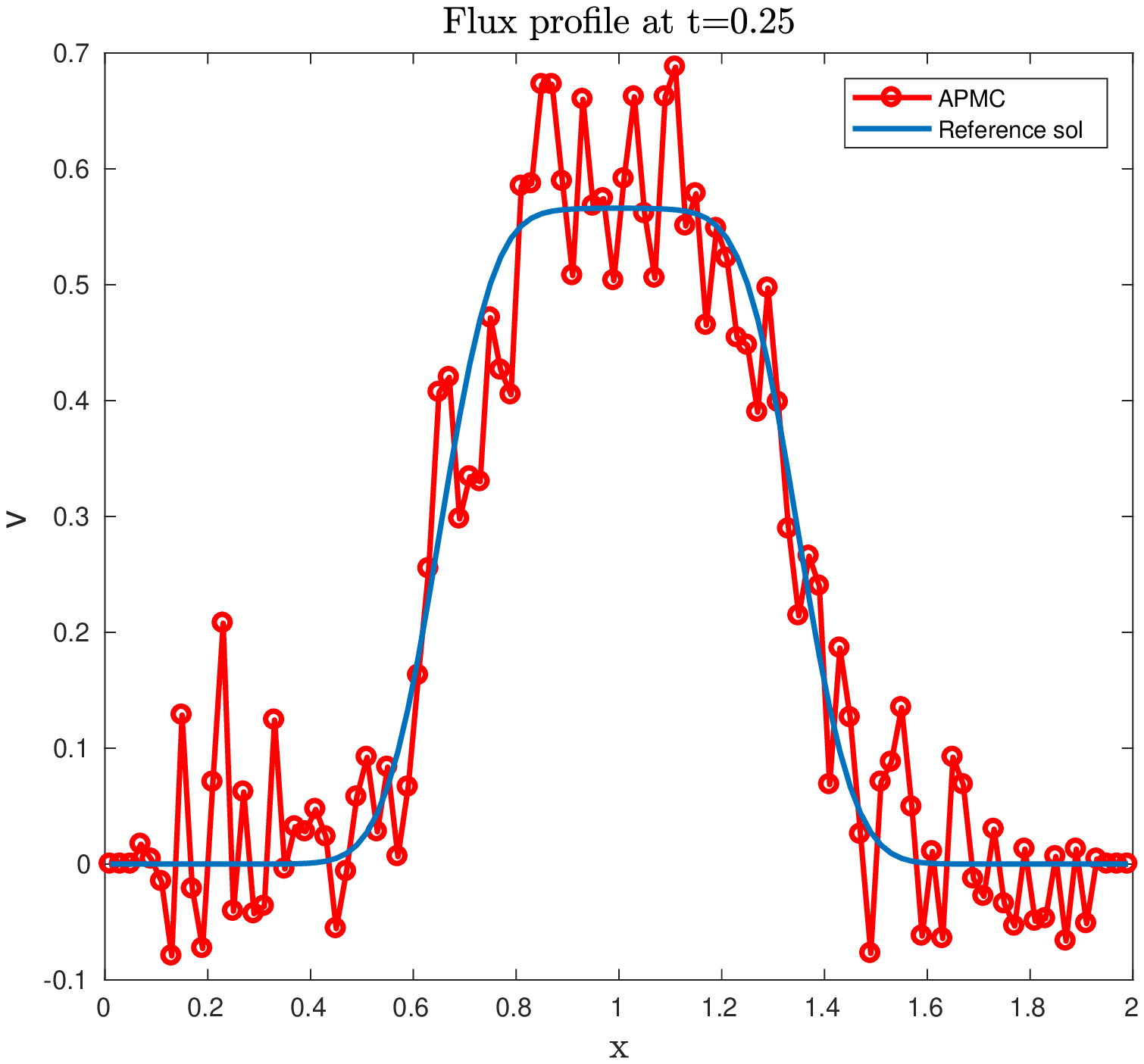}\hspace{1cm}
          \includegraphics[scale=0.37]{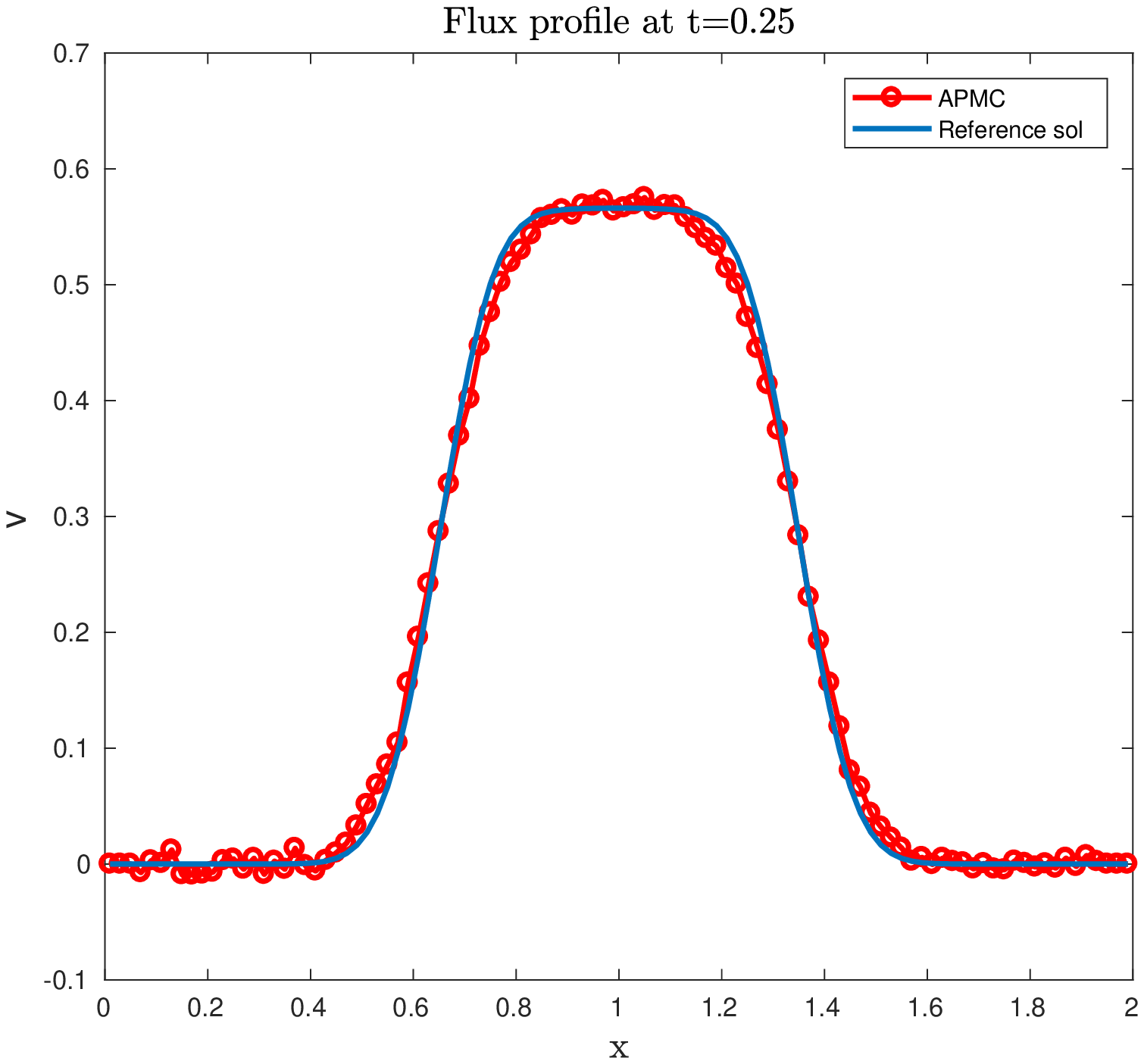}
\caption{Goldstein-Taylor. Numerical solution at time t = 0.25 in the rarefied regime $\varepsilon= 0.7$
with $\Delta t = 0.01$ and $\Delta x = 0.02$. The mass density $u$ (top) and the flow $v$ (bottom) are shown (red circles) together with a reference solution (blue continuous line). On the left panels
the solution of the Monte Carlo scheme with $1000$ particles per cell in average is shown, on the right panels the converged solution is represented.
}
\label{fig1}
\end{center}
\end{figure}

\begin{figure}[ht!]
\begin{center}
    \includegraphics[scale=0.37]{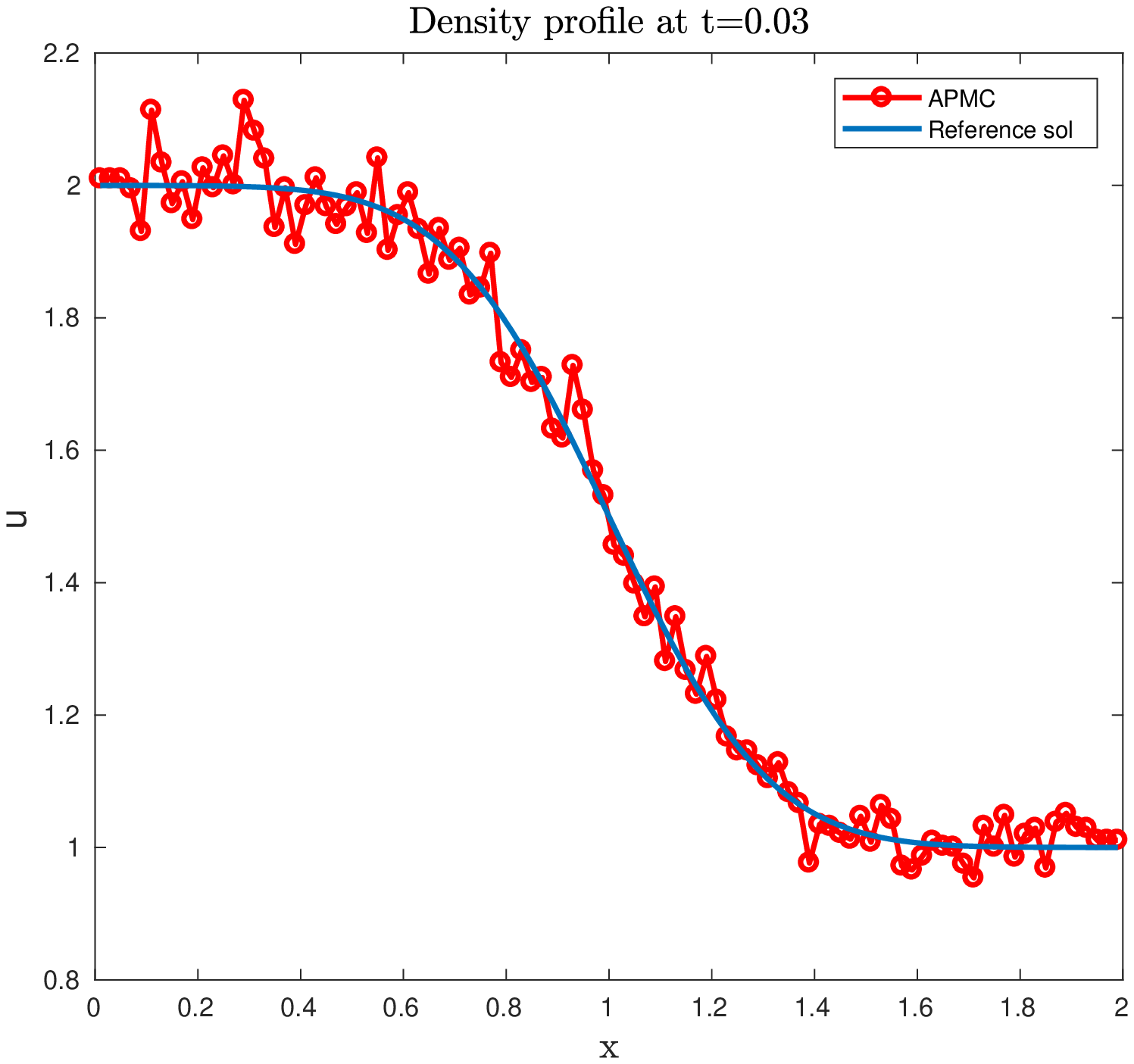}\hspace{1cm}
        \includegraphics[scale=0.37]{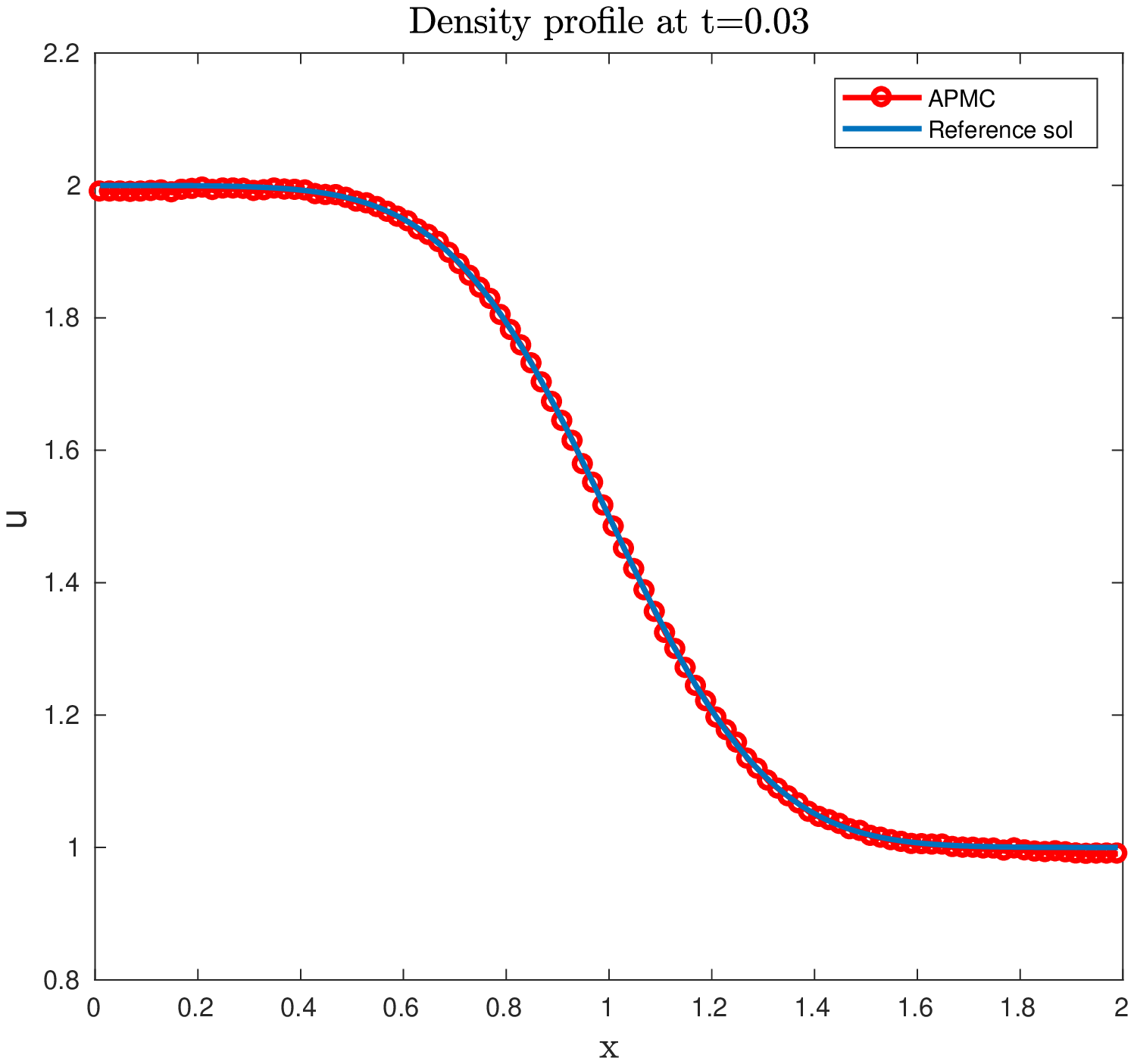}\\
     \includegraphics[scale=0.37]{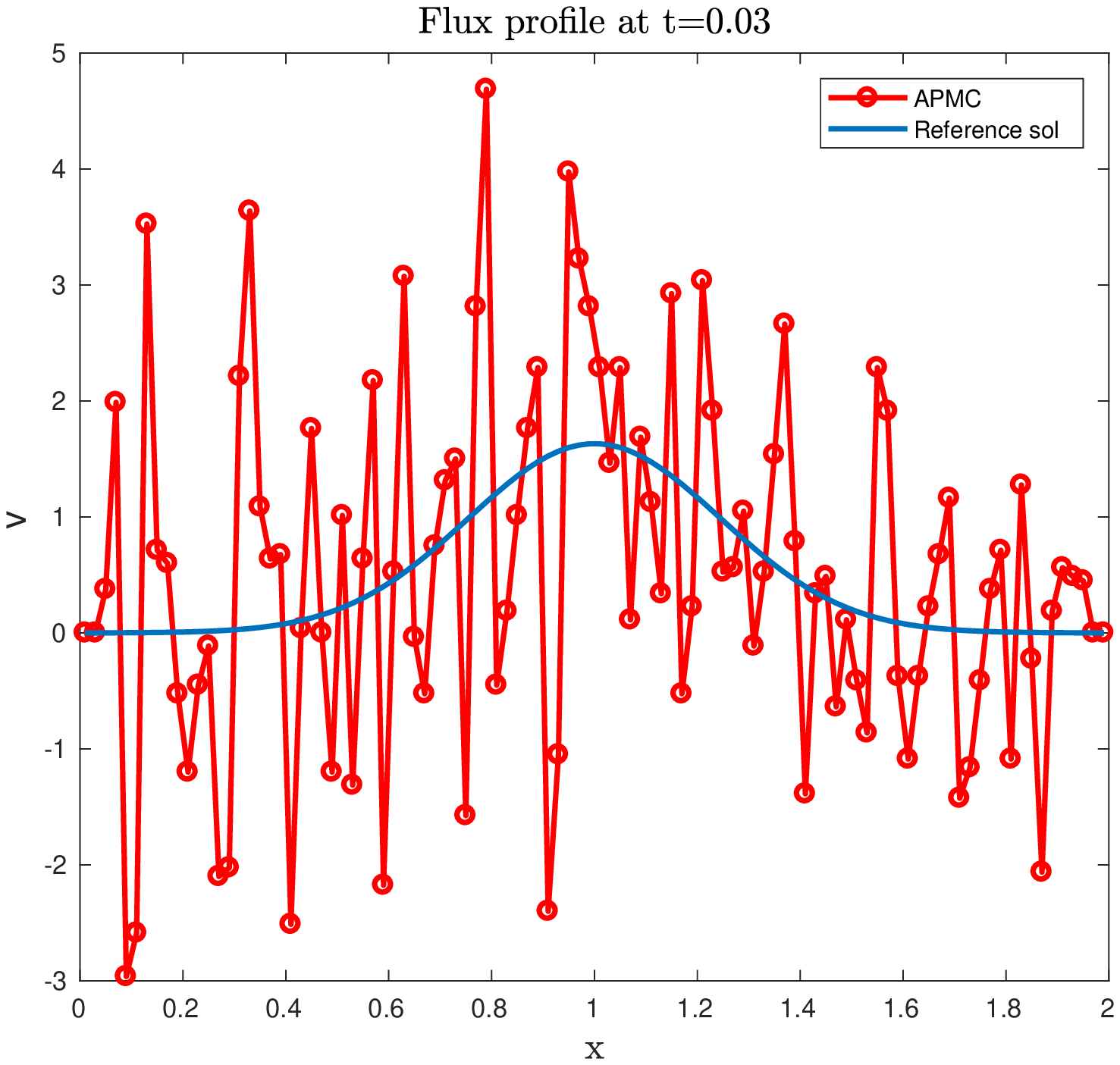}\hspace{1cm}
          \includegraphics[scale=0.37]{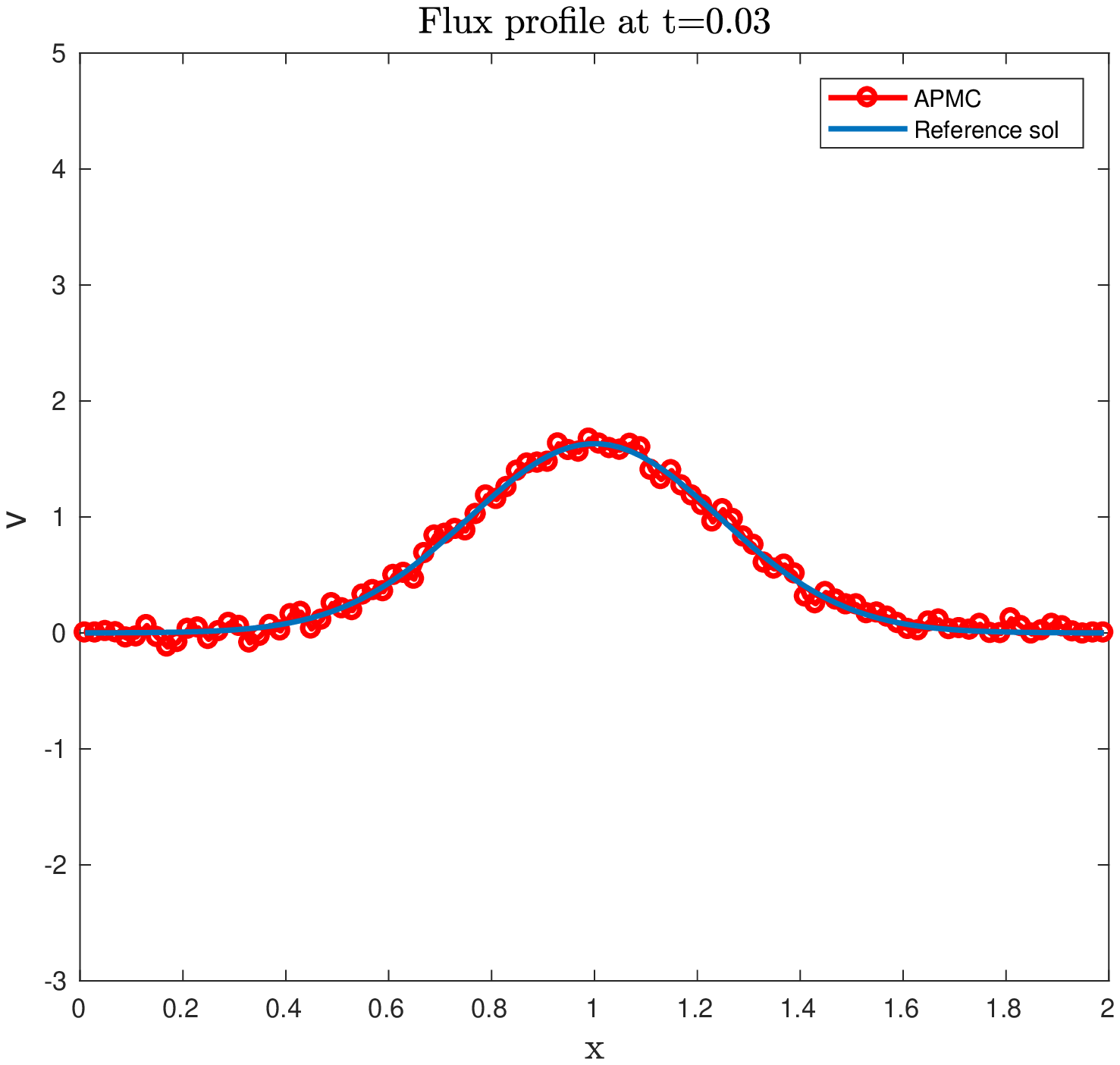}
\caption{Goldstein-Taylor. Numerical solution at time t = 0.03 in the rarefied regime $\varepsilon= 10^{-5}$
with $\Delta t = 0.01$ and $\Delta x = 0.02$. The mass density $u$ (top) and the flow $v$ (bottom) are shown (red circles) together with a reference solution (blue continuous line). On the left panels
the solution of the Monte Carlo scheme with $1000$ particles per cell in average is shown, on the right panels the converged solution is represented.
}
\label{fig2}
\end{center}
\end{figure}

We consider the following Riemann problem with initial data 
\[
\rho_L = 2.0, \qquad j_L=0 \qquad 0<x<1,\]
\[\rho_R = 1.0, \qquad j_R=0 \qquad 1<x<2.
\]
The test is run both in the diffusive, i.e. $\varepsilon=10^{-5}$, as well as in the hyperbolic regime, i.e. $\varepsilon=0.7$. This test models the behaviors of to two semi-infinite rods having different initial temperature and put into contact at initial time. The boundary conditions are of Dirichlet type and fix the density $\rho_L(x=0,t)$ and $\rho_R(x=2,t)$ and the flux $j_L(x=0,t)$ and $j_R(x=2,t)$ for all times. An analogous problem has been studied in \cite{JPT1}. In Figure (\ref{fig1}) and (\ref{fig2}) the solution obtained with the Monte Carlo method with $N=100$ mesh points is reported together with a reference solution which employs the same number of mesh points. In the simulations shown, the time step is for $\varepsilon=10^{-5}$, $\Delta t=0.4 \ (\Delta x)^2$ while for $\varepsilon=0.7$ is $\Delta t=0.5\ \Delta x$. Top images show the density profiles, while bottom images show the flux profiles. Left images show the solution computed with an average of $1000$ particles per cell while right images show the converged solutions in terms of the number of particles. The convergence of the density function is much faster than the convergence of the fluxes. This is due to the definition of density and flux functions: the first is defined as sum of positive and negative mass particles while the second is defined as the difference of the number of positive and negative particles divided by the scaling parameter $\varepsilon$. This means that when $\varepsilon$ becomes very small the convergence of the flux function becomes very hard. This can be observed in the Figure \ref{fig2} on the bottom right for which still some fluctuations are present even if the density is fully converged. 
 In both cases, the numerical solutions match the reference solution very well. 

\subsection{The radiative transport}\label{numerics-RT}
\subsubsection{One dimensional case}
\begin{figure}[ht!]
\begin{center}
    \includegraphics[scale=0.37]{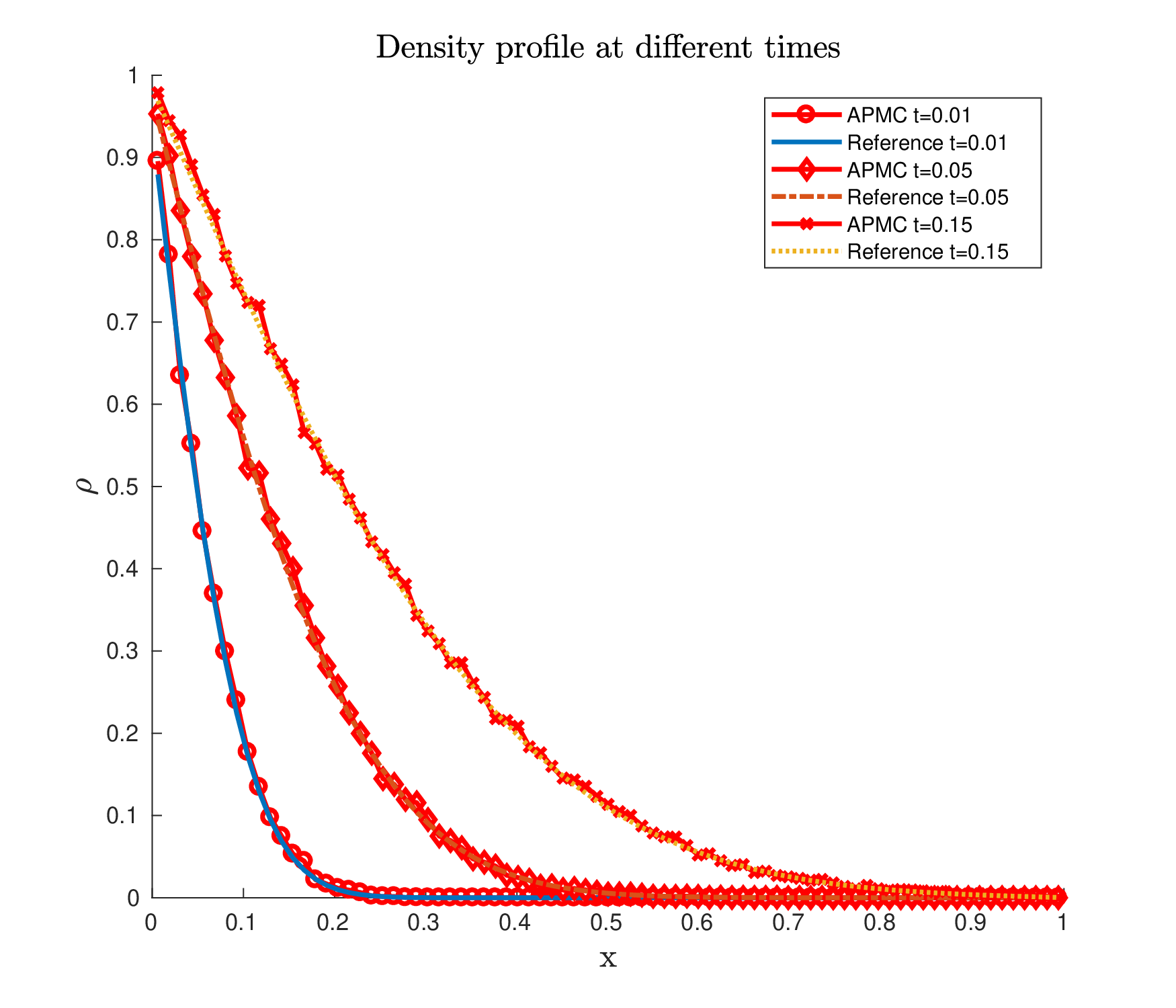}\hspace{1cm}
        \includegraphics[scale=0.37]{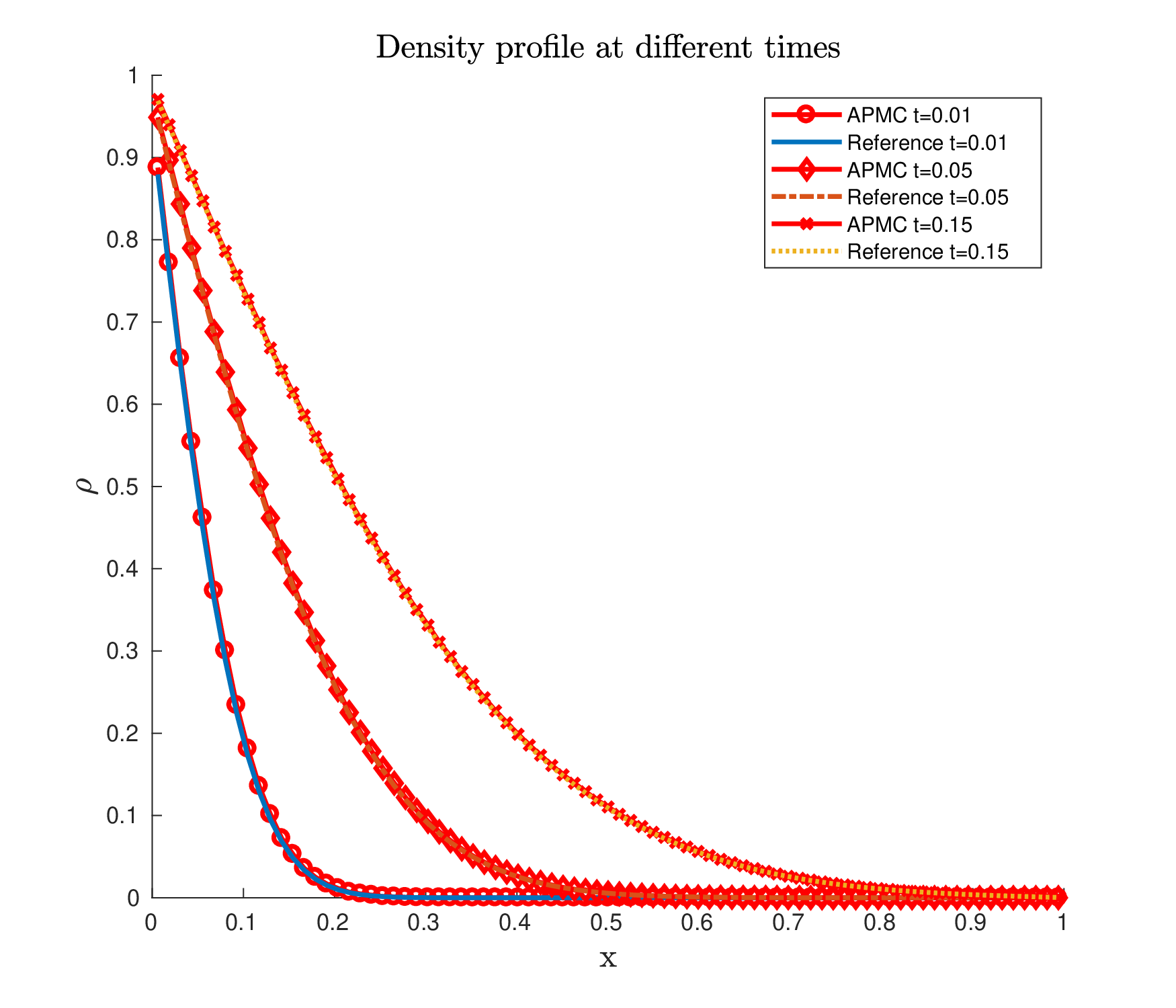}\\
        \includegraphics[scale=0.37]{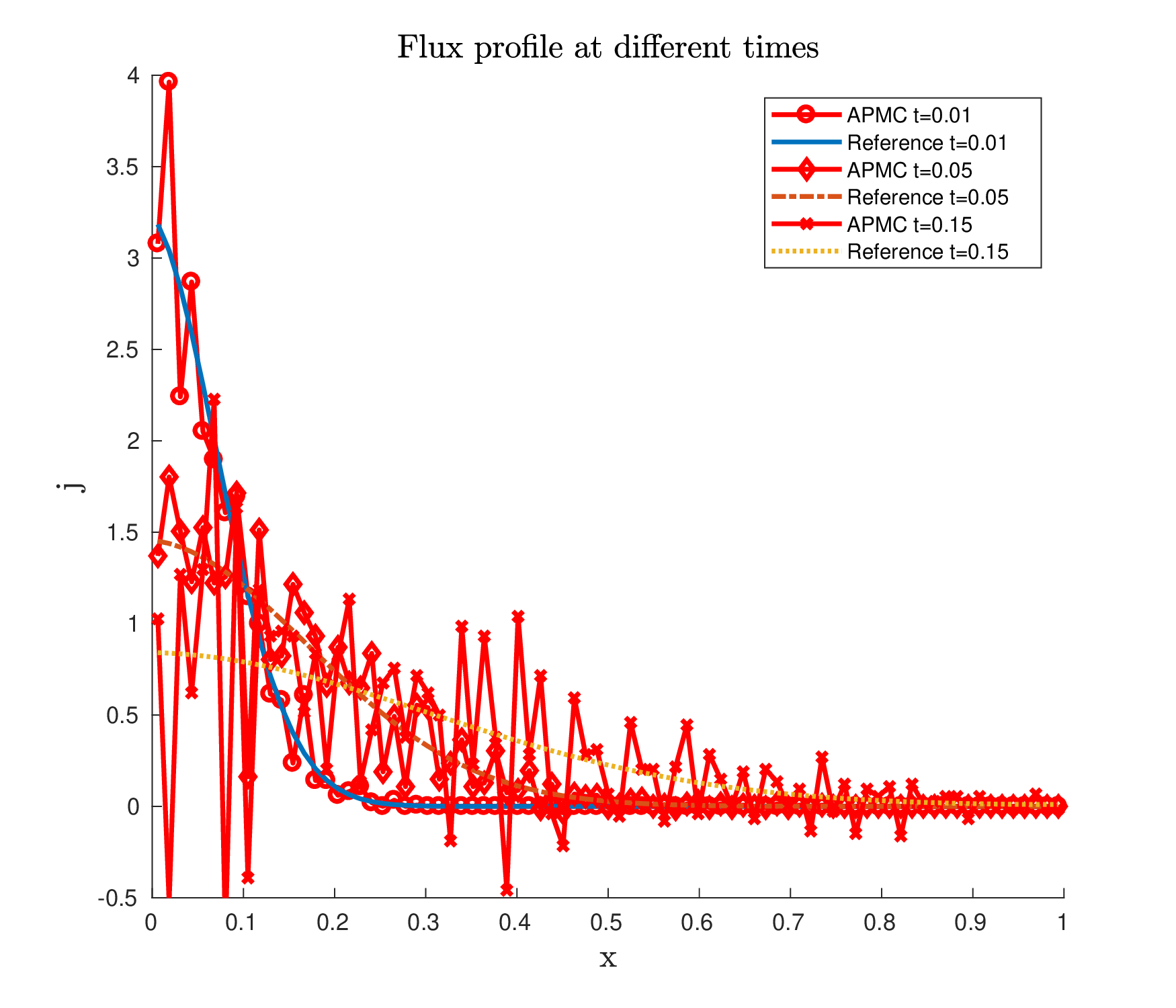}\hspace{1cm}
        \includegraphics[scale=0.37]{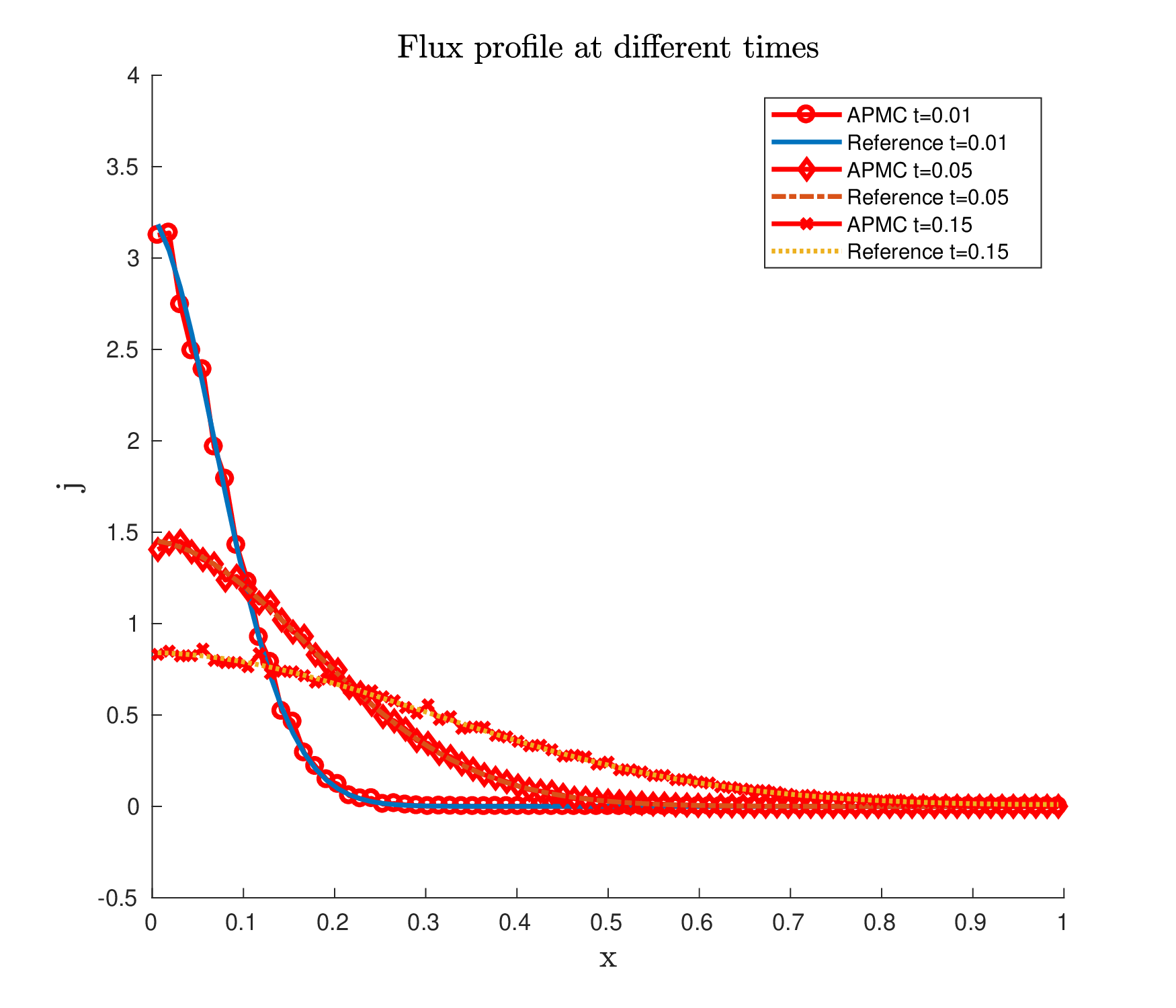}
\caption{Radiative transport (Test 1). Numerical solutions at time $t = 0.01$, $t=0.05$ and $t=0.15$ in the diffusive regime $\varepsilon= 10^{-8}$
with $\Delta t = 5 \ 10^{-4}$ and $\Delta x = 0.01$. The mass density $\rho$ (top) and the flow $j$ (bottom) are shown. Red circles for $t=0.01$, red diamonds for $t=0.05$ and red squares for $t=0.15$.  Blue continuous, dash dotted and dotted lines for reference solutions. Left pictures report the solution obtained with $1000$ particles in average per cell while right pictures report the converged solutions. }
\label{fig4}
\end{center}
\end{figure}

\begin{figure}[ht!]
\begin{center}
    \includegraphics[scale=0.5]{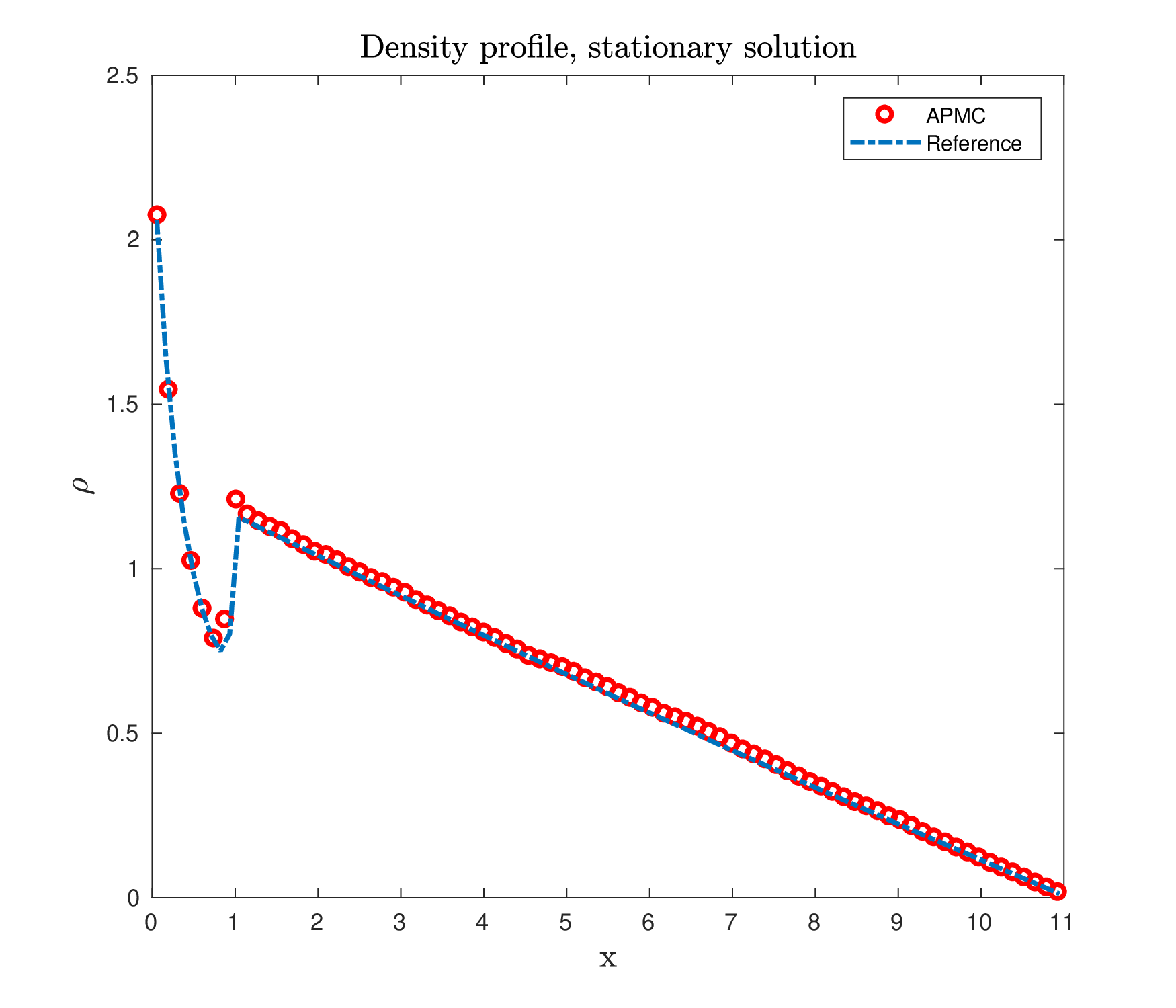}
\caption{Radiative transport (Test 2). Steady state numerical solution $\varepsilon= 1$ in the rarefied region, $\varepsilon=0.01$ in the dense region.
Red circle for the Monte Carlo solution, blue continuous line reference solution.
}
\label{fig3}
\end{center}
\end{figure}

We consider two transport problems in slab geometry. The first problem has the following initial and boundary data
\[x \in [0, 1], \ \sigma_S = 1, \ \varepsilon=10^{-8} ,\ \sigma_A = 0,\] \[ f(x=1,v,t)= 0,\ f(x=0,v,t)= 1, \ f(x,v,t)=0 \ \forall x\in(0,1).
\]
The test is run in the diffusive regime, i.e. $\varepsilon=10^{-8}$. An analogous problem has been studied for instance in \cite{JPT2}. In Figure (\ref{fig4}), the solution obtained with the Monte Carlo method with $N=80$ mesh points is reported together with a reference solution which employs the same number of points. The images show the solution at different instant of time, namely $t=0.01$, $t=0.05$ and $t=0.15$. Top images show the density profiles, while bottom images show the flux profiles. The left images report a solution obtained with $1000$ particles in average per cell, while the right images report the converged solution. As for the two speed case, the flux function is measured by the difference of positive and negative particles speed divided by the scaling factor $\varepsilon$. This means that in the pure diffusive regimes these are difficult to obtain by means of particle schemes due to the very fine resolution demanded. The time step is fixed to $\Delta t=0.5\ (\Delta x)^2$. 
The scheme is able to furnish correct solutions even for choices of time and space step which are much larger than the scaling parameter $\varepsilon$.
The second problem considered has the following initial and boundary data \cite{JPT2}
\[ x \in [0, 11], \ f(x=0,v,t)=5, \ f(x=11,v,t)=0 \]
\[ \sigma_S = 0, \ \sigma_A=1, \ \varepsilon= 1, \ x \in ]0, 1],\]
\[ \sigma_S = 1, \ \sigma_A=0, \ \varepsilon= 0.01, \ x \in [1, 11].\]
In the purely absorbing region the solution decays exponentially whereas in the purely
scattering region the solution is diffusive. The solution is computed by using $N=80$ mesh points both for the Monte Carlo and the reference solution and it is run up to a stationary solution is reached. For the intermediate states $100$ particles per cell are used in average. The time step is the one of the diffusive region, i.e. $\Delta t=0.5 (\Delta x)^2$.
Time average techniques has been used to reduce the statistical noise by averaging successive solutions once the steady state is reached. At the interface between the two regions, for the particles which cross the interface, the transport and the diffusion coefficients are the ones on the left of the interface before the particle reach this point and the ones on the right after the interface has been reached. The numerical scheme gives a good description of the internal layer and of the solution in the absorption and diffusive regions as shown in Figure \ref{fig3}. 

\subsubsection{Two dimensional case}
\begin{figure}[ht!]
\begin{center}
    \includegraphics[scale=0.34]{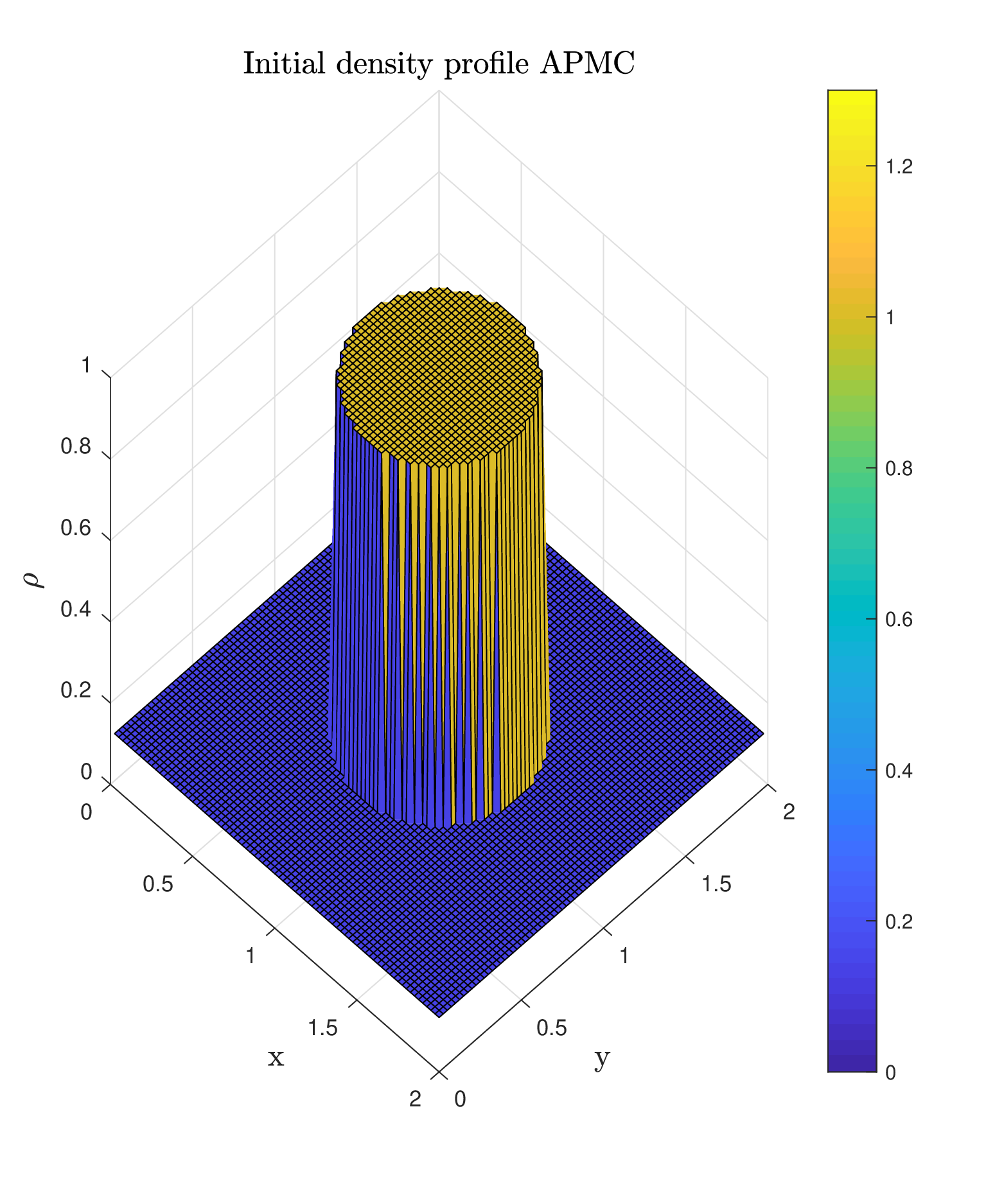}\hspace{1cm}
        \includegraphics[scale=0.34]{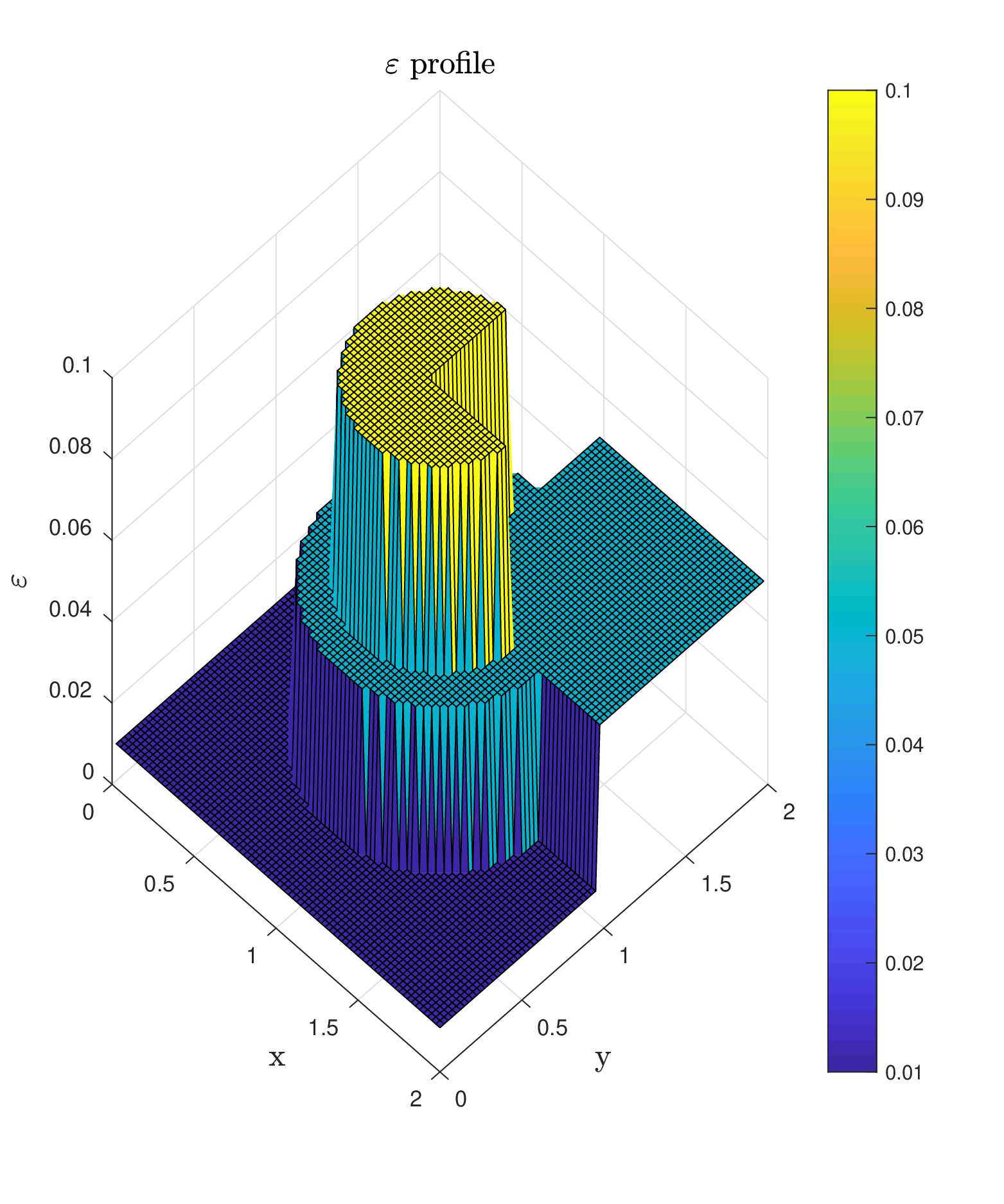}\\
             \includegraphics[scale=0.34]{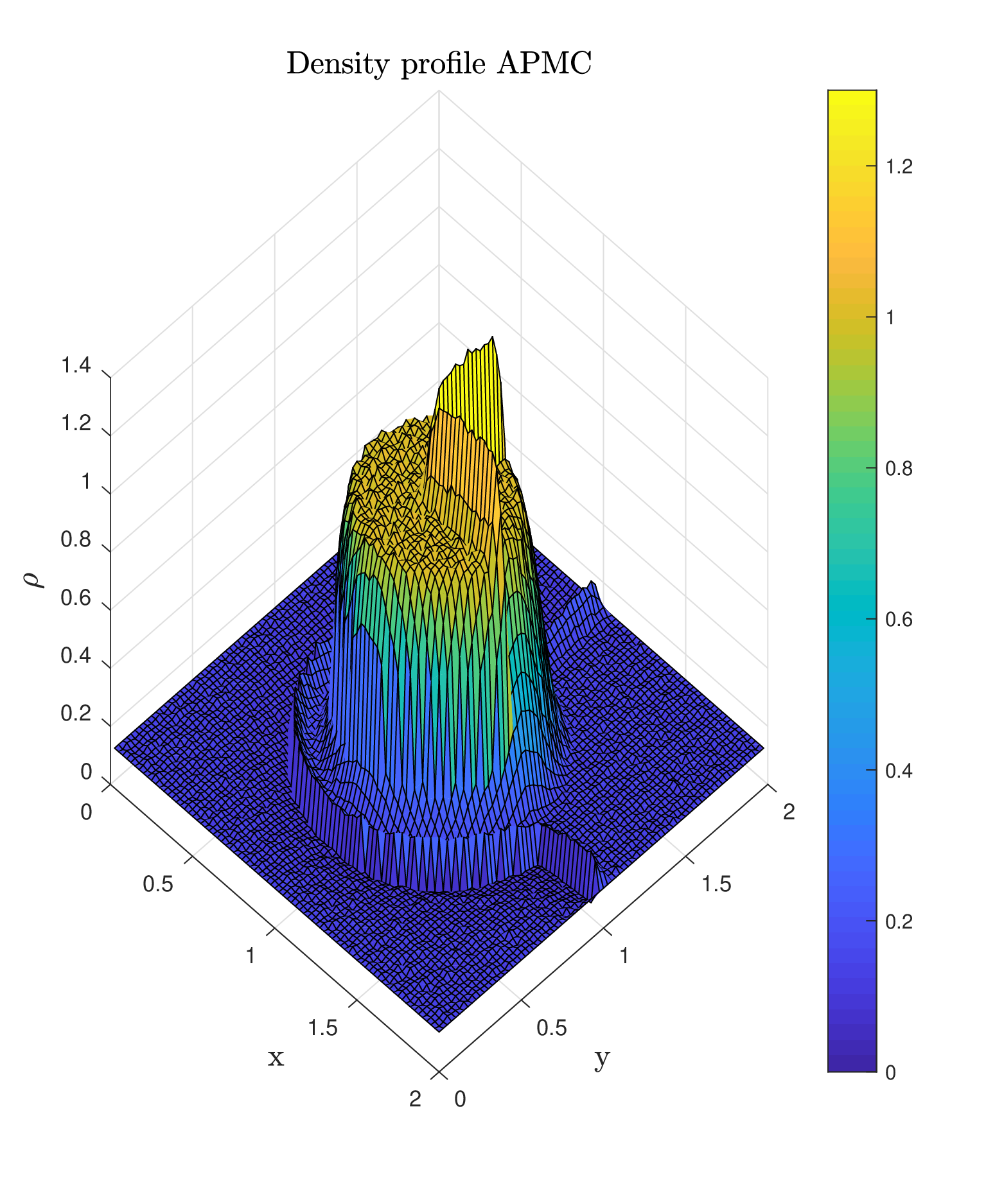}\hspace{1cm}
          \includegraphics[scale=0.34]{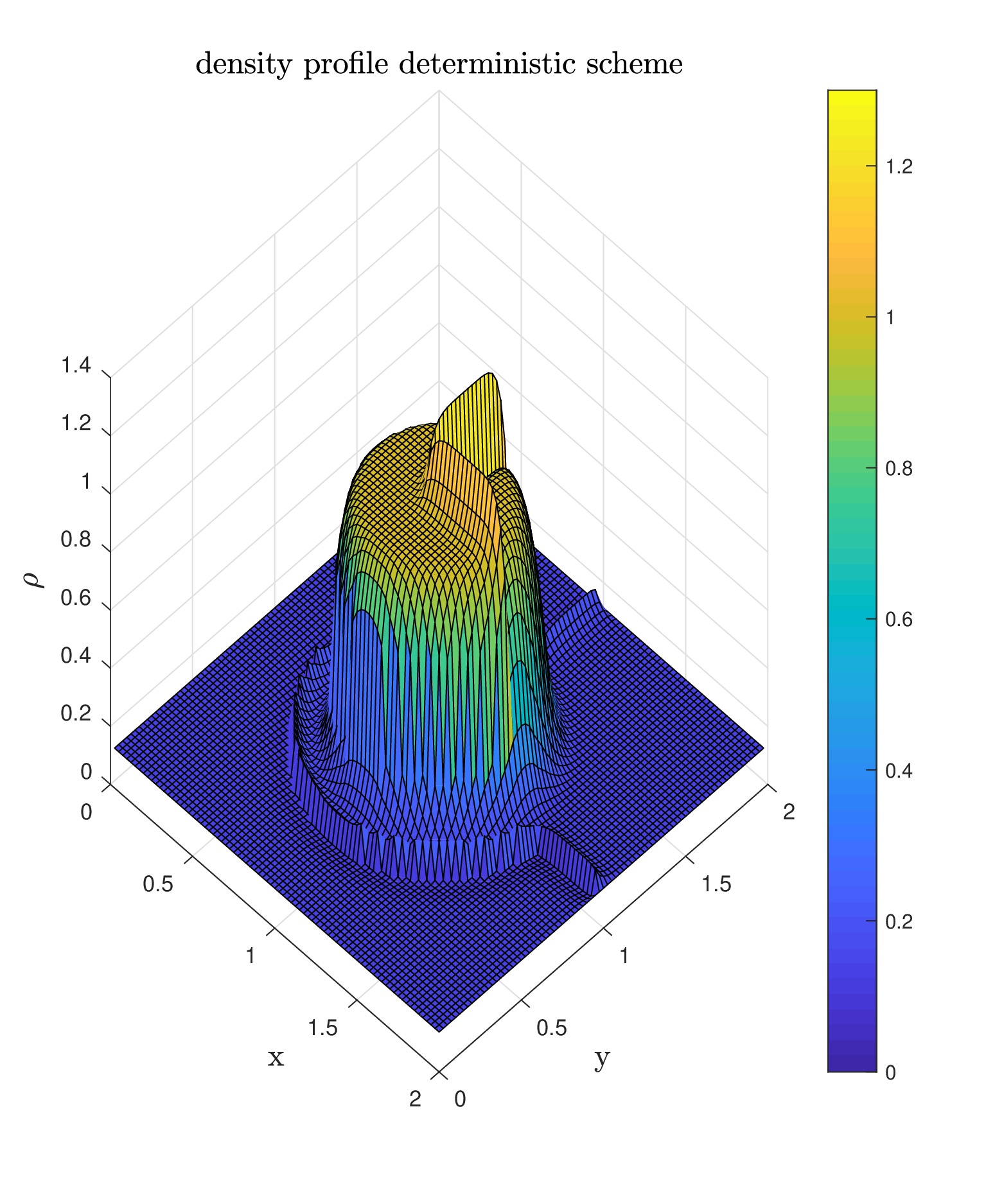}
\caption{Radiative transport in two dimensions. Initial density profile top left, $\varepsilon$ profile top right.
Numerical solution at time t = 0.025 for the asymptotic preserving Monte Carlo method bottom left, Numerical solution at time t = 0.025 for a reference solution bottom left.
The Monte Carlo scheme employs $2000$ particles per cell in average.
}
\label{fig5}
\end{center}
\end{figure}

We consider one transport problem in two dimensions in space and velocity space in bounded domain $[0,2]^2$. As initial data we consider an uniform distribution function $f$ such that $\int f dv=1$ in a central circular region of radius $0.2$, while $\int f dv=0.125$ outside of this region. In Figure \ref{fig5} on the top left such initial data is shown. The transport coefficients are set to $\sigma_s=1$ everywhere without absorption, i.e. $\sigma_a=0$. The scaling coefficient $\varepsilon$ is discontinuous ranging from $0.1$ to $0.01$ and it is represented in Figure \ref{fig5} on the top right. The number of cells in space are $80\times 80$ while the time step is $\Delta t=\Delta x^2$ and the final time of the simulation is fixed to $T=0.002$. The boundary conditions are of Dirichlet type. The solution obtained with the Monte Carlo method is reported together with a deterministic  solution which employs the same number of points in space and $20\times 20$ in velocity space in Figure \ref{fig5} respectively on the bottom left and right. The number of particles is an average $2000$ per cell. The images show that the MC method proposed is able to well describe a varying relaxation regime without a time step dependent on the scaling parameter $\varepsilon^2$.

\section{Discussion and conclusions}\label{Concl}
A new class of Monte Carlo schemes for solving transport equations in the diffusive limit has been presented. The approach is based on a reformulation of the original equations
in order to obtain a modified system which characteristic speeds do not arbitrary grow when the scaling parameter goes to zero. The idea is to introduce a suitable implicit time discretization for the original model that permits to derive an equivalent system (up to a first order error in time) with bounded characteristic speeds. The resulting Monte Carlo schemes are unconditionally stable with respect the scaling parameter and degenerate automatically in the limit to a classical Brownian Monte Carlo solver for the diffusive equation without any time step limitations. In the last part, several numerical tests have been performed which show the capability of the method to deal with different situations from rarefied to diffusive regimes. In the next future extension of this methodology to other diffusion limits, like semiconductor kinetic equations, and the construction of hybrid schemes which combine the Monte Carlo solver with a deterministic solver for the limiting diffusion equation will be considered.

\section*{Acknowledgements}

This research was made possible by the Royal Flemish Academy of Belgium for Sciences and Art (KVAB), which supported a visit of GS to the University of Ferrara via a travel grant. We also acknowledge the support of the INDAM-GNCS grant ''Uncertainty quantification for hyperbolic and kinetic equations''.

\bibliographystyle{plain}

\end{document}